\documentclass[a4paper,9pt]{article}

\usepackage[centertags]{amsmath}

\usepackage{authblk}
\usepackage{multirow}
\usepackage{amssymb}
\usepackage{graphicx}
\usepackage{listing}
\usepackage{mathrsfs}
\usepackage{amsmath,amsthm,amssymb,fontenc,color,textcomp,amsfonts,graphicx,graphics}
\usepackage[colorlinks=true]{hyperref}
\usepackage{hyperref}
\hypersetup{colorlinks = true,linkcolor = black,citecolor = blue}
\usepackage{xcolor}
\usepackage{booktabs}

\theoremstyle{definition}
\newtheorem{defn}{Definition}[section]
\newtheorem{thm}{Theorem}[section]

\newtheorem{Ex}{Example}[section]

\usepackage[linesnumbered,ruled,vlined]{algorithm2e}

\usepackage{caption}
\usepackage{subcaption}

\newcommand{\be}{\begin{equation}}
\newcommand{\ee}{\end{equation}}

\begin{document}
\date{}
\title{Tensor regularization by truncated iteration: a comparison of some solution methods for large-scale linear discrete ill-posed problem with a t-product\\}

\author{{\Large Ugochukwu O. Ugwu\thanks{\,e-mail: uugwu@kent.edu} \;and\; Lothar Reichel\thanks{\,e-mail: reichel@math.kent.edu}}
\mbox{\ }\\
{\normalsize \large Department of Mathematical Sciences, Kent State University, OH 44242, USA}}


\maketitle \vspace*{-0.5cm}

\thispagestyle{empty}

\begin{abstract}
\noindent
This paper describes and compares some structure preserving techniques for the solution of linear discrete ill-posed problems with the t-product. A new randomized tensor singular value decomposition (R-tSVD) with a t-product is presented for low tubal rank tensor approximations. Regularization of linear inverse problems by truncated tensor eigenvalue decomposition (T-tEVD), truncated tSVD (T-tSVD), randomized T-tSVD (RT-tSVD), t-product Golub-Kahan bidiagonalization (tGKB) process, and t-product Lanczos (t-Lanczos) process are considered. A solution method that is based on reusing tensor Krylov subspaces generated by the tGKB process is described. The regularization parameter is the number of iterations required by each method. The discrepancy principle is used to determine this parameter. Solution methods that are based on truncated iterations are compared with solution methods that combine Tikhonov regularization with the tGKB and t-Lanczos processes. Computed examples illustrate the performance of these methods when applied to image and gray-scale video restorations. Our new RT-tSVD method is seen to require less CPU time and yields restorations of higher quality than the T-tSVD method. \vspace{.3cm}

\noindent
{\bf key words:} linear discrete ill-posed problem, tensor Golub-Kahan bidiagonalization, t-singular triplets, t-eigen pairs, randomized tensor eigenvalue decomposition,  tensor singular value decomposition, tensor Lanczos process.
\end{abstract}

\section{Introduction}
\noindent
The solution of linear systems
\[
\mathcal{A}*\mathcal{\vec{X}} = \mathcal{\vec{B}},
\]
has been recently investigated in literature; see, e.g., \cite{KBHH, KM, RU1, RU2, RU3, RU4}, with significant research efforts devoted to the solution of large-scale least squares problems,
\begin{equation}
\min_{\mathcal{\vec{X}} \in \mathbb{R}^{m \times 1 \times n}} \|\mathcal{\vec{B}} - \mathcal{A}*\mathcal{\vec{X}}\|_F, ~~~ \mathcal{A} \in \mathbb{R}^{\ell \times m \times n}, ~~~ \mathcal{\vec{B}} \in \mathbb{R}^{\ell \times 1 \times n}, ~~~ \ell \geq m.
\label{I3}
\end{equation}

We are concerned with the solution of \eqref{I3} when the third order tensor $\mathcal{A} = [a_{ijk}]_{i,j,k=1}^{\ell,m,n}$ is of ill-determined tubal rank, and the Frobenius norm of the singular tubes of $\mathcal{A}$ decay rapidly to zero with increasing index; see, e.g., \cite{RU1, RU2}. The singular tubes of $\mathcal{A}$ are analogues of the singular values of a matrix, and there are many nonvanishing singular tubes of tiny Frobenius norm of different orders of magnitude close to zero. Then \eqref{I3} is referred to as a linear discrete ill-posed problem. The tensor slices $\mathcal{\vec{X}}$ and $\mathcal{\vec{B}}$ are laterally oriented $m \times n$ and $\ell \times n$ matrices, respectively; see Section \ref{sec2}. The operator $*$ described below denotes the t-product introduced in the seminal work by Kilmer and Martin \cite{KM}. This product has become ubiquitous in tensor literature applications; see, e.g., facial recognition \cite{HKBH}, tomographic image reconstruction \cite{SKH}, video completion \cite{ZEAH}, image classification \cite{NKH}, and image deblurring \cite{GIJS, KBHH, KM, RU1, RU2, RU3, RU4}. Throughout this paper, $\|\mathcal{A} \|_F$ denotes the Frobenius norm of a third order tensor $\mathcal{A}$. It is given by
\[
\|\mathcal{A}\|_F= \sqrt{\sum_{i=1}^\ell \sum_{j=1}^{m} \sum_{k=1}^n a^2_{ijk}}.
\]

The data tensor slice $\mathcal{\vec{B}} \in \mathbb{R}^{\ell \times 1 \times n}$ in our applications, e.g., image and video restorations, is contaminated by measurement error represented by the tensor slice, $\mathcal{\vec{E}} \in \mathbb{R}^{\ell \times 1 \times n}$, such that
\begin{equation*}
\mathcal{\vec{B}} = \mathcal{\vec{B}}_\text{true}+\mathcal{\vec{E}},
\end{equation*}
where $\mathcal{\vec{B}}_\text{true}\in \mathbb{R}^{\ell \times 1 \times n}$ denotes the unknown error free data tensor slice. Assume that the unavailable system of equations
\[
\mathcal{A}*\mathcal{\vec{X}}= \mathcal{\vec{B}}_\text{true}
\]
is consistent with an exact solution, $\mathcal{\vec{X}}_\text{true} \in \mathbb{R}^{m \times 1 \times n}$. A direct solve for  the solution of \eqref{I3} does not yield a meaningful approximation of $\mathcal{\vec{X}}_\text{true}$ due to severe ill-conditioning of $\mathcal{A}$, and error propagation into the solution of \eqref{I3}. However, we can recover a meaningful approximate solution by replacing \eqref{I3} with a nearby problem whose solution is less sensitive to the error in $\mathcal{\vec{B}}$. This replacement is known as regularization.

A common regularization technique in matrix theory is the truncated singular value decomposition (TSVD). The tensor analogue of the TSVD for regularizing \eqref{I3} is the truncated tensor singular value decomposition (T-tSVD) introduced by Kilmer and Martin \cite{KM}. This approach truncates the tSVD of $\mathcal{A}$ to $k$-terms and yields the best tubal rank-$k$ approximation $\mathcal{A}_k$ of $\mathcal{A}$. The optimality of the truncated tSVD is established in \cite{KM}. The truncation index $k$ is the regularization parameter determined by the discrepancy principle. Assume that the bound
\begin{equation} \label{errbd}
\|\mathcal{\vec{E}}\|_F \leq \delta,
\end{equation}
is available. Let $\mathcal{\vec{X}}_k$ denote an approximate solution of \eqref{I3}. The discrepancy principle prescribes that the smallest integer $k > 0$ be determined such that
\begin{equation} \label{disp}
\|\mathcal{\vec{B}} - \mathcal{A}*\mathcal{\vec{X}}_k \|_F \leq \tau \delta,
\end{equation}
where $\tau > 1$ is a user specified constant independent of $\delta>0$. For details of the discrepancy principle as well as other methods that can be applied to determine the regularization parameter $k$; see, e.g., \cite{EHN, FRR,Ha}. 

It is the purpose of this paper to discuss and compare some structure preserving techniques for the solution of linear discrete ill-posed problems with the t-product. An advantage of preserving tensor structures when solving \eqref{I3} is to avoid information loss inherent due to flattening; see \cite{KBHH}. Computed examples presented in  \cite{RU1,RU2,RU3} illustrate the merits of tensorizing over matricizing or vectorizing ill-posed tensor equations. Solution methods described by El Guide et al. \cite{GIJS} involve flattening, i.e., they reduce \eqref{I3} to an equivalent equation involving a matrix and a vector. Structure preserving and other techniques for regularizing \eqref{I3} by Tikhonov's approach with the t-product are described and compared in \cite{RU1, RU2, RU3, RU4}. We will discuss regularization of \eqref{I3} by truncated iterations with the tensor eigenvalue decomposition (tEVD), tSVD, randomized tSVD (R-tSVD), t-product Golub-Kahan bidiagonalization (tGKB) process, and t-product Lanczos (t-Lanczos) process.

We describe how a few of the tensor singular (t-singular) triplets that are made up of the singular tubes of the largest Frobenius norm and their corresponding left and right singular matrices can be computed inexpensively. Linear discrete ill-posed problems \eqref{I3} when $\mathcal{A}$ is a small matrix $A$ and $\mathcal{\vec{B}}$ is a vector $b$ are commonly solved by computing the SVD of $A$; see, e.g., \cite{Ha,NR}. However, for a general large-scale matrix $A$, computing the SVD is expensive. Thus, one focuses on computing a few of the largest singular triplets of $A$. These triplets are made up of the largest singular values and the associated left and right singular vectors. It is shown recently in \cite{OR} that the computation of a few of the largest singular triplets of $A$ by the implicitly restarted Golub-Kahan bidiagonalization process generally is quite inexpensive, and much cheaper than computing the (full) SVD of $A$. Similarly, for a symmetric matrix $A$, the eigenpairs of the largest magnitude that consist of a few eigenvalues of the largest magnitude and the corresponding eigenvectors can be computed less fairly inexpensively by the implicitly restarted Lanczos process; see, e.g., \cite{BCR1, BCR2,B,OR}.

Inspired by the results presented in \cite{KBHH,KM,OR}, we describe how a few of the largest t-singular triplets of $\mathcal{A}$ can be computed inexpensively by using the tGKB process introduced in \cite{KBHH}. The tGKB process generates an orthonormal tensor basis for the t-product Krylov (t-Krylov) subspace
\begin{equation}\label{krylov1}
\mathbb{K}_k(\mathcal{A}^T*\mathcal{A},\mathcal{A}^T*\mathcal{\vec{B}})=
{\rm span}\{\mathcal{A}^T*\mathcal{\vec{B}},
(\mathcal{A}^T*\mathcal{A})*\mathcal{A}^T*\mathcal{\vec{B}},\ldots,
(\mathcal{A}^T*\mathcal{A})^{k-1}*\mathcal{A}^T*\mathcal{\vec{B}}\},
\end{equation}
where $^T$ denotes transposition. This process requires two tensor-matrix product evaluations with $\mathcal{A}$ and with $\mathcal{A}^T$, and reduces the tensor $\mathcal{A}$ in \eqref{I3} to a small $(k+1) \times k \times n$ lower bidiagonal tensor by carrying out a few, $k \ll m$, steps. 

When a symmetric tensor $\mathcal{A} \in \mathbb{R}^{m \times m \times n}$ defines the problem \eqref{I3}, the t-Lanczos process described in \cite{KBHH} can be applied to express the solution of \eqref{I3} as a tensor linear (t-linear) combination of the lateral slices of a Krylov tensor generated by $\mathcal{A}$ and $\mathcal{\vec{B}}$ in the t-Krylov subspace
\begin{equation}\label{krylov2}
\mathbb{K}_\ell(\mathcal{A},\mathcal{\vec{B}})=
{\rm span}\{\mathcal{\vec{B}},
\mathcal{A}*\mathcal{\vec{B}},\ldots,
\mathcal{A}^{\ell-1}*\mathcal{\vec{B}}\}.
\end{equation}
The t-Lanczos process generates an orthonormal tensor basis for the t-Krylov subspace \eqref{krylov2}, and computes a few $\ell \ll m$ steps to reduce $\mathcal{A}$ to a small $(\ell +1) \times \ell \times n$ tridiagonal tensor. This process may be advantageous over the tGKB process in terms of computational cost and storage since only a tensor-matrix product evaluation with $\mathcal{A}$ is required.

We describe how the t-Lanczos process applied to $\mathcal{A}$ is used to compute a few
eigentubes of the largest Frobenius norm and associated eigenmatrices which we refer to as the tensor eigenpairs (t-eigenpairs) of the largest magnitude. The tensor eigenvalue decomposition (tEVD) of
$\mathcal{A} \in \mathbb{R}^{m \times m \times n}$ has been described in \cite{KBHH,LJ}
using the t-product. The solution of \eqref{I3} by truncated tEVD (T-tEVD) of $\mathcal{A}$ is discussed in
Subsection \ref{sec3.1}.

We also consider solution methods for the minimization problem analogous to \eqref{I3},
\begin{equation}\label{I6}
\min_{\mathcal{X} \in \mathbb{R}^{m \times p \times n}} \|\mathcal{B} - \mathcal{A}*\mathcal{X}\|_F, ~~~ \mathcal{A} \in \mathbb{R}^{\ell \times m \times n}, ~~~ \mathcal{B} \in \mathbb{R}^{\ell \times p \times n}, ~~~ \ell \geq m,~~~ p>1,
\end{equation}
where the tensors $\mathcal{X}$ and $\mathcal{B}$ are general third order tensors. These kinds of problems arise in color image and video restorations. Solution methods for \eqref{I6} by Tikhonov regularization with the t-product Arnoldi (t-Arnoldi) and tGKB processes are discussed in \cite{GIJS, RU1,RU2,RU3}. 

We consider the regularization of \eqref{I6} by the T-tEVD, T-tSVD, t-Lanczos and tGKB methods. The solution method for \eqref{I6} based on nested t-Krylov subspaces with the tGKB process described in \cite{RU1} is also considered. The data tensor $\mathcal{B}$ has lateral slices $\mathcal{\vec{B}}_j$, $j=1,2,\dots, p$. We will work with the lateral slices of $\mathcal{B}$ independently to determine an approximate solution of \eqref{I6}. The resulting solution methods for \eqref{I6} based on this approach are referred to as the tGKB$_p$ and t-Lanczos$_p$ methods. Both methods apply the tGKB and t-Lanczos methods $p$ times to solve \eqref{I6}. 

Randomized tSVD algorithms introduced by Zhang et al. \cite{ZSKA} can also be applied to determined
an approximate solution of \eqref{I3} and \eqref{I6}. The algorithms in \cite{ZSKA} are well suited for tubal
rank problems of known low tubal rank, i.e., when we have prior knowledge of the tubal rank of $\mathcal{A}$. Since this is often not the
case when solving linear discrete ill-posed problems, we address this difficulty by presenting a new fixed-precision-based
R-tSVD algorithm that relies on a set accuracy tolerance to determine a low tubal rank approximation of $\mathcal{A}$.
The proposed R-tSVD algorithm is based on the randomized block algorithms for matrix computations presented in
\cite{MV, YGL}. The randomized solution of linear discrete ill-posed problems defined by a matrix and a vector has
received some attention in literature; see, eg., \cite{BBR,JY,XZ1,XZ2}. For a complete survey of randomized
solution methods for large-scale problems in two dimensions; see \cite{HMT}.

The structure of this paper is as follows. Section 2 introduces notations and preliminaries associated with the t-product, while Section \ref{sec3} discusses the solution of \eqref{I3} and \eqref{I6} when $\mathcal{A}$ is symmetric. The computation of a few t-eigenpairs of $\mathcal{A}$ using the t-Lanczos process is discussed in Subsection \ref{sec3.2}. In Section \ref{sec4}, solution methods for \eqref{I3} and \eqref{I6} are considered for a nonsymmetric tensor $\mathcal{A}$. An approach that is based on the T-tSVD is outlined in Section \ref{sec4.1}. Low tubal rank approximation of $\mathcal{A}$ by the R-tSVD is discussed in Subsection \ref{sec4.2}. Derivation of the solutions of \eqref{I3} and \eqref{I6} in terms of the largest t-singular triplets of $\mathcal{A}$, by using the tGKB process, are presented in Subsection \ref{sec4.3}. Also considered in this subsection is the nested t-Krylov subspace method that is based on the tGKB process. 

Computed examples illustrate the performance of the described methods in Section \ref{sec5} when applied to image and gray-scale video restorations.  Comparisons between solution methods that are based on the tGKB and t-Lanczos processes with truncated iterations and Tikhonov regularization are presented. The Tikhonov's approach requires more computational time and often computes solution of better quality than methods based on truncated iterations.  Randomization can be used to speed up and improve the quality 
of the computed solution determined by the T-tSVD method. The proposed randomized T-tSVD (RT-tSVD) method requires less CPU time and yields restorations of higher quality than the T-tSVD method. Section \ref{sec6} gives the concluding remarks.

\section{Notations and Preliminaries}\label{sec2}

A tensor is a multidimensional array of numbers. In this paper, a tensor is of order three and will be denoted by calligraphic script, say, $\mathcal{A} \in \mathbb{R}^{\ell \times m \times n}$. The order is the number of modes or ways of the tensor. Matrices and vectors are tensors of second and first order, respectively. We will use capital letters, say, $A$ to denote matrices. For third order tensors, a slice is a $2$D section obtained by fixing one of the indices. Using MATLAB notation, $\mathcal{A}(i,:,:)$, $\mathcal{A}(:,j,:)$ and $\mathcal{A}(:,:,k)$ denotes the $i$th horizontal, $j$th lateral, and $k$th frontal slices, respectively. The $j$th lateral slice denoted by $\mathcal{\vec{A}}_j$ is a tensor referred to as the tensor column, whereas the $k$th frontal slice $\mathcal{A}^{(k)}$ is a matrix. A fiber of a third order tensor is a $1$D section obtained by fixing two of the indices. Using MATLAB notation, $\mathcal{A}(:,i,j)$, $\mathcal{A}(i,:,j)$, and $\mathcal{A}(i,j,:)$ denote the $(i,j)$th mode-1, mode-2 and mode-3 fibers, respectively. We denote the tubal fibers (tubes or mode-3 fibers) by a bold-face letter, say, ${\bf a}$. Definition \ref{tprod} gives the definition of the t-product.

\begin{defn}{(t-product \cite{KM})}\label{tprod}
Suppose $\mathcal{B} \in \mathbb{R}^{\ell \times p \times n}$ and $\mathcal{C} \in \mathbb{R}^{p \times m \times n}$. Then the t-product between the tensors $\mathcal{B}$ and $\mathcal{C}$ is the $\ell \times m \times n$ tensor $\mathcal{A}$ whose $(i,j)$th tube is given by
\begin{equation*}
\mathcal{A}(i,j,:) = \sum_{k=1}^p \mathcal{B}(i,k,:)\ast \mathcal{C}(k,j,:),
\end{equation*}
where $*$ denotes the circular convolution between tubes.
\end{defn}
The t-product uses the Discrete Fourier Transform (DFT) to efficiently compute the circular convolution between tubes. Throughout this paper, we use ~$\widehat{\cdot}$~ to indicate objects obtained by taking FFT along the tubes of a third order tensor. Let the Fast Fourier Transform (FFT) along the tubes of $\mathcal{B}$ and $\mathcal{C}$ be $\mathcal{\widehat{B}} = \mathtt{fft}(\mathcal{B},[\;],3)$ and $\mathcal{\widehat{C}} = \mathtt{fft}(\mathcal{C},[\;],3)$, respectively. The t-product is computed by performing a matrix-matrix product of each pair of the frontal slices of $\mathcal{\widehat{B}}$ and $\mathcal{\widehat{C}}$ in the Fourier domain, i.e.,
\begin{equation*}
\mathcal{\widehat{A}}(:,:,i) = \mathcal{\widehat{B}}(:,:,i) \cdot \mathcal{\widehat{C}}(:,:,i), ~~~~ i = 1, 2, \dots, n,
\end{equation*}
and taking the inverse FFT of the result along the third dimension to obtain $\mathcal{A} = \mathtt{ifft}(\mathcal{\widehat{A}},[\;],3)$. It can be shown that taking FFT of $\mathcal{A} \in \mathbb{R}^{\ell \times m \times n}$ along the third dimension will cost $\mathcal{O}(\ell m n {\rm log}_2(n))$; see \cite{KM}.

The tensor transpose has similar property as the matrix transpose. Given $\mathcal{A} \in \mathbb{R}^{\ell \times m \times n}$, then $\mathcal{A}^T \in \mathbb{R}^{m \times \ell \times n}$ is the tensor obtained by transposing each of the frontal slices and then reversing the order of the transposed frontal slices 2 through $n$; see \cite{KM}.
It can be shown that if $\mathcal{A}$ and $\mathcal{B}$ are two tensors such that $\mathcal{A*B}$ and $\mathcal{B}^T*\mathcal{A}^T$ are defined. Then $(\mathcal{A*B})^T = \mathcal{B}^T*\mathcal{A}^T$; see \cite{KM}.

The identity tensor $\mathcal{I} \in \mathbb{R}^{m \times m \times n}$ is a tensor whose first frontal slice, $\mathcal{I}^{(1)}$ is the $m \times m$ identity matrix and all other frontal slices, $\mathcal{I}^{(i)}$, $i=2,3,\dots,n$, are zero matrices; see \cite{KM}. Related to the identity tensor is the tensor $\vec{e}_1 \in \mathbb{R}^{m \times 1 \times n}$ whose $(1,1,1)$th entry is equals $1$ and the remaining entries are zeros. The nonzero entry $1$ appear only at the first frontal face of $\vec{e}_1$. The tensor ${\bf e}_1 \in \mathbb{R}^{1 \times 1 \times n}$ is tube with $(1,1,1)$th entry equals 1 and the remaining entries varnish; see \cite{ZA1}.

A tensor $\mathcal{Q} \in \mathbb{R}^{m \times m \times n}$ is orthogonal if $\mathcal{Q}^T*\mathcal{Q} = \mathcal{Q}*\mathcal{Q}^T = \mathcal{I}$; see \cite{KM}. Analogous to the columns of orthogonal matrices, the lateral slices of $\mathcal{Q}$ are orthonormal. The concept of partial orthogonality is similar to a tall thin matrix that has orthonormal columns. In the case of a tensor, $\mathcal{Q} \in \mathbb{R}^{\ell \times m \times n}$, $\ell > m$, is partially orthogonal if $\mathcal{Q}^T*\mathcal{Q}$ is well defined and equal to the identity tensor $\mathcal{I}_{mmn}$; see \cite{KM}. The lateral slices of $\mathcal{Q}$
are orthonormal, i.e.,
\begin{equation*}
\mathcal{Q}^T(:,i,:)*\mathcal{Q}(:,j,:) = \left\{
\begin{array}{ll}
{\bf e}_1 & i=j,\\
{\bf 0} & i \neq j,
\end{array}
\right.
\end{equation*}
where $\mathbf{e}_1\in\mathbb{R}^{1\times 1\times n}$ is a tubal scalar with
$(1,1,1)$th entry equal to $1$ and the remaining entries are zero. It can be shown that if $\mathcal{Q}$ is an orthogonal tensor, then
\be \label{normF}
\| \mathcal{Q*A}\|_F = \|\mathcal{A}\|_F;
\ee
see \cite{KM}. An $m \times m \times n$ tensor $\mathcal{A}$ has an inverse denoted by $\mathcal{A}^{-1}$, provided that $\mathcal{A}*\mathcal{A}^{-1} = \mathcal{I}$ and $\mathcal{A}^{-1}*\mathcal{A} = \mathcal{I}$. A tensor is called f-diagonal if each frontal slice is a diagonal matrix; see \cite{KM}.

Given $\mathcal{A} \in \mathbb{R}^{\ell \times m \times n}$. The tSVD of $\mathcal{A}$ is given by
\begin{equation*}
\mathcal{A} = \mathcal{U*S}*\mathcal{V}^T,
\end{equation*}
where $\mathcal{U} \in \mathbb{R}^{\ell \times \ell \times n}$ and $\mathcal{V} \in \mathbb{R}^{m \times m \times n}$ are orthogonal tensors, and $\mathcal{S} \in \mathbb{R}^{\ell \times m \times n}$ is an f-diagonal tensor; see \cite{KM}. Two important measures of tensor complexity from the tSVD that are used in dimensionality reduction are the tensor multi-rank and tubal rank. The multi-rank of $\mathcal{A} \in \mathbb{R}^{\ell \times m \times n}$ is a tube {\bf k}, whose $i$th face corresponds to the rank of the $i$th frontal slice of $\mathcal{A}$. Thus, the multi-rank of $\mathcal{A}$ is $k$ if ${\bf k}^{(i)} = k$, $i = 1,2,\dots,n$; see \cite{KBHH}. Let $\mathcal{A} = \mathcal{U*S}*\mathcal{V}^T$. Then the tensor tubal rank denoted by $\mathtt{rank}_t(\mathcal{A})$, is the number of nonzero singular tubes of $\mathcal{S}$, i.e.,
\begin{equation*}
\mathtt{rank}_t(\mathcal{A}) = \#\{i: \mathcal{S}(i,i,:) \neq 0\} = \max\limits_{i} {\bf k}^{(i)}.
\end{equation*}

We say $\mathcal{A} \in \mathbb{R}^{\ell \times m \times n}$ has full tubal rank if $\mathtt{rank}_t(\mathcal{A}) = \min\{\ell,m\}$. Two interesting properties of tensor tubal rank that are similar to matrix rank are
\[
\mathtt{rank}_t(\mathcal{A})\leq \min(\ell,m) ~~{\rm and}~~ \mathtt{rank}_t(\mathcal{A}*\mathcal{B})\leq \min(\mathtt{rank}_t(\mathcal{A}),\mathtt{rank}_t(\mathcal{B}));
\]
see \cite{JN}. Suppose $\mathcal{A} \in \mathbb{R}^{\ell \times m \times n}$ has tubal rank-$k$. Then the truncated tSVD of $\mathcal{A}$ described by Algorithm \ref{Alg: tsvd} yields the factorization
\begin{equation*}
\mathcal{A}_k = \mathcal{U}_k*\mathcal{S}_k*\mathcal{V}_k^T,
\end{equation*}
where $\mathcal{U}_k \in \mathbb{R}^{\ell \times k \times n}$ and $\mathcal{V}_k \in \mathbb{R}^{m \times k \times n}$ are partially orthogonal tensors, and $\mathcal{S} \in \mathbb{R}^{k \times k \times n}$ is f-diagonal; see \cite{KM}.

\vspace{.3cm}
\begin{algorithm}[H]
\SetAlgoLined
\KwIn{$\mathbf{Input}: \mathcal{A} \in \mathbb{R}^{\ell \times m \times n}$}
\KwOut{$\mathcal{U} \in \mathbb{R}^{\ell \times k \times n},\;\; \mathcal{S} \in \mathbb{R}^{k \times k \times n},\;\; \mathcal{V} \in \mathbb{R}^{m \times k \times n}$}
$\mathcal{\widehat{A}} = \mathtt{fft}(\mathcal{A},[\;],3)$\\
\For {$i = 1$ \bf{to} $n$}{
$[U_k, S_k, V_k] = \mathtt{svd}\mathtt(\mathcal{\widehat{A}}^{(i)},k)$\\
$\mathcal{\widehat{U}}_k^{(i)} = U_k$, $\mathcal{\widehat{S}}_k^{(i)} = S_k$, $\mathcal{\widehat{V}}_k^{(i)} = V_k$\\
}
$ \mathcal{U}_k= \mathtt{ifft}(\mathcal{\widehat{U}}_k,[\;],3),\; \mathcal{S}_k= \mathtt{ifft}(\mathcal{\widehat{S}}_k,[\;],3),\; \mathcal{V}_k= \mathtt{ifft}(\mathcal{\widehat{V}}_k,[\;],3)$
\caption{Truncated tensor singular value decomposition (T-tSVD) \cite{KM}.}
\label{Alg: tsvd}
\end{algorithm}\vspace{.3cm}

Let $\mathcal{A}\in \mathbb{R}^{\ell \times m \times n}$. The following theorem shown in \cite{KM} discusses the optimality of the error in the truncated tSVD of $\mathcal{A}$.

\begin{thm}
Let the tSVD of $\mathcal{A}\in \mathbb{R}^{\ell \times m \times n}$ be given by $\mathcal{A} = \mathcal{U*S}*\mathcal{V}^T$. For $k \leq \min\{\ell, m\}$, define
\begin{equation*}
\mathcal{A}_k = \sum_{j=1}^k \mathcal{\vec{U}}_j*\mathbf{s}_j*\mathcal{\vec{V}}_j^T,
\end{equation*}
where $\mathcal{\vec{U}}_j = \mathcal{U}(:,j,:)$ and $\mathcal{\vec{V}}_j = \mathcal{V}(:,j,:)$ are the lateral slices of $\mathcal{U}$ and $\mathcal{V}$, respectively, and $\mathbf{s}_j = \mathcal{S}(j,j,:)$ is the $(j,j)$th tube of $\mathcal{S}$. Then $\mathcal{A}_k = \arg\min\limits_{\mathcal{\tilde{A}}\in \mathcal{M}} \|\mathcal{A}-\mathcal{\tilde{A}}\|_F$, where $\mathcal{M} = \{\mathcal{X*Y} | \mathcal{X}\in \mathbb{R}^{\ell \times k \times n}, \mathcal{Y}\in \mathbb{R}^{k \times m \times n}\}$. Moreover, $\| \mathcal{A} - \mathcal{A}_k\|_F$ is the minimal error over all tensors of tubal rank-$k$. It is given by
\[
\| \mathcal{A} - \mathcal{A}_k\|_F = \Bigg(\frac{1}{n} \sum_{i=1}^n \sum_{j=k+1}^{\min\{\ell,m\}} \widehat{\sigma}_j^{(i)} \Bigg)^{1/2},
\]
where $\widehat{\sigma}_j^{(i)} = \widehat{\mathcal{S}}(j,j,i)$ is the $j$th singular value corresponding to the $i$th face of $\widehat{\mathcal{A}}$.
\label{thm 2.1}
\end{thm}
We see from Theorem \ref{thm 2.1} that $\mathcal{A}_k$ is a tensor of tubal rank at most $k$, and is a best tubal rank-$k$ approximation to $\mathcal{A}$. We conclude this subsection with two interesting concepts associated with the tSVD; the range (tensor column space) of a tensor $\mathcal{A}$, and the orthogonal projection onto the range of $\mathcal{A}$.

Suppose the tSVD factorization $\mathcal{A} = \mathcal{U*S}*\mathcal{V}^T$ has tubal rank-$k$. The range of $\mathcal{A}$, denoted by $\mathcal{R}(\mathcal{A})$, is the space spanned by the lateral slices $\mathcal{\vec{U}}_j$, $j=1,2,\dots,k$, of $\mathcal{U}$. The set
\begin{equation*}
\mathcal{R}(\mathcal{A}) = \{ \mathcal{\vec{U}}_1*\mathbf{c}_1 + \dots + \mathcal{\vec{U}}_k*\mathbf{c}_k\; | \;\mathbf{c}_j \in \mathbb{R}^{1 \times 1 \times n} \}, ~~~j = 1,2,\dots,k,
\end{equation*}
is determined by the t-linear combination of the lateral slices of $\mathcal{U}$ that form an orthonormal tensor basis for the range of $\mathcal{A}$; see \cite{NKH}. Given a lateral slice $\mathcal{\vec{B}} \in \mathbb{R}^{\ell \times 1 \times n}$, the orthogonal projection onto the range of $\mathcal{A}$ is defined by $ \mathcal{U}*\mathcal{U}^T*\mathcal{\vec{B}}$; see \cite{NKH}.

\section{Symmetric t-linear discrete ill-posed problems}\label{sec3}
This section reviews the tensor eigenvalue decomposition (tEVD) and tensor Lanczos process described by Kilmer et al. \cite{KBHH} using the t-product, and discusses how they can be applied to determine an approximate solution of \eqref{I3} and \eqref{I6}.

\subsection{The tEVD method for the approximate solution of \eqref{I3} and \eqref{I6}}\label{sec3.1}

Let $\mathcal{A} \in \mathbb{R}^{m\times m\times n}$ in \eqref{I3} be symmetric with tEVD \cite{KBHH},
\begin{equation}
\mathcal{A} = \mathcal{W}* \mathcal{D} *\mathcal{W}^T \implies \mathcal{A}*\mathcal{\vec{W}}_j =
\mathcal{\vec{W}}_j*{\mathbf{d}}_j, ~~ j = 1, 2,\ldots, m,
\label{eq: 3}
\end{equation}
where $\mathcal{W} = [\mathcal{\vec{W}}_1, \mathcal{\vec{W}}_2, \dots, \mathcal{\vec{W}}_m] \in
\mathbb{R}^{m\times m\times n}$ is an orthogonal tensor with $m$ orthonormal tensor columns, and
$\mathcal{D} = \text{diag}[{\mathbf{d}}_1, {\mathbf{d}}_2, \dots, {\mathbf{d}}_m]
\in \mathbb{R}^{m\times m\times n}$ is an f-diagonal tensor with eigentubes ${\mathbf{d}}_j$
corresponding to the eigenmatrices $\mathcal{\vec{W}}_j$. Assume that the eigentubes ${\mathbf{d}}_j$
decay rapidly and are ordered\footnote{\url{https://www.mathworks.com/matlabcentral/fileexchange/18904-sort-eigenvectors-eigenvalues}} according to
\[
\|{\mathbf{d}}_1\|_F \geq \|{\mathbf{d}}_2\|_F \geq \dots \geq \|{\mathbf{d}}_m\|_F.
\]
Then only the first few eigentubes of the largest Frobenius norm and associated eigenmatrices are of interest. We remark that  $\mathcal{A}$ is symmetric if each frontal slice
$\mathcal{\widehat{A}}^{(i)}$, $i = 1,2,\dots, n$, of $\mathcal{\widehat{A}}$ is Hermitian; see \cite{KBHH}.

Define the truncated tEVD computed by Algorithm \ref{tevd} as
\begin{equation}
\mathcal{A}_k = \mathcal{W}_k* \mathcal{D}_k *\mathcal{W}^T_k,
\label{eq: 4}
\end{equation}
where $\mathcal{W}_k= [\mathcal{\vec{W}}_1, \mathcal{\vec{W}}_2, \dots, \mathcal{\vec{W}}_k]
\in \mathbb{R}^{m\times k\times n}$ and $\mathcal{D}_k = \text{diag}[{\mathbf{d}}_1, {\mathbf{d}}_2,
\dots, {\mathbf{d}}_k] \in \mathbb{R}^{k\times k\times n}$ for some $1 \leq k \leq m$. By Theorem
\ref{thm 2.1}, $\mathcal{A}_k$ is the closest tensor to $\mathcal{A}$ in the Frobenius norm among all
tensors of tubal rank at most $k$. Thus, $\mathcal{A}_k$ is the best tubal rank-$k$ approximation of
$\mathcal{A}$ computed by Algorithm \ref{tevd}. The eigenvalues and corresponding eigenvectors of
$\mathcal{\widehat{A}}^{(i)}$, $i=1,2,\dots,n$, of $\mathcal{\widehat{A}}$ are sorted in a descending
order of magnitude in step $3$ of Algorithm \ref{tevd}.

\vspace{.3cm}
\begin{algorithm}[H]
\SetAlgoLined
\KwIn{$\mathbf{Input}: \mathcal{A} \in \mathbb{R}^{m \times m \times n}$}
\KwOut{$\mathcal{W} \in \mathbb{R}^{m \times k \times n},\;\; \mathcal{D} \in \mathbb{R}^{k \times k \times n}$}
$\mathcal{\widehat{A}} = \mathtt{fft}(\mathcal{A},[\;],3)$\\
\For {$i = 1$ \bf{to} $n$}{
$[W_k, D_k] = {\tt sort}\big(\mathtt{eig}\mathtt(\mathcal{\widehat{A}}^{(i)},k),{\rm 'descend'}\big)$\\
$\mathcal{\widehat{W}}_k^{(i)} = W_k$, $\mathcal{\widehat{D}}_k^{(i)} = D_k$\\
}
$ \mathcal{W}_k= \mathtt{ifft}(\mathcal{\widehat{W}}_k,[\;],3),\; \mathcal{D}_k= \mathtt{ifft}(\mathcal{\widehat{D}}_k,[\;],3)$
\caption{Truncated tensor eigenvalue decomposition (T-tEVD)}
\label{tevd}
\end{algorithm}\vspace{.3cm}

Replacing $\mathcal{A}$ with $\mathcal{A}_k$, and $\mathcal{\vec{B}}$ with $\mathcal{W}_k* \mathcal{W}_k^T * \mathcal{\vec{B}}$ (the orthogonal projector onto the range of $\mathcal{A}_k$) in \eqref{I3} for a suitable (small) value of $k$ determined by the discrepancy principle \eqref{disp}, and letting $\mathcal{\vec{X}} = \mathcal{W}_k * \mathcal{\vec{Y}} \in \mathbb{R}^{k \times 1\times n}$, we obtain the reduced minimization problem
\begin{equation}
\min_{\mathcal{\vec{Y}}} ||\mathcal{D}_k * \mathcal{\vec{Y}} - \mathcal{W}^T_k * \mathcal{\vec{B}}||_F.
\label{eq: 5}
\end{equation}
The solution of \eqref{eq: 5} is given by $\mathcal{\vec{Y}}_k = \mathcal{D}^{-1}_k*\mathcal{W}^T_k * \mathcal{\vec{B}}$, and the approximate solution of \eqref{I3} can be expressed as
\begin{equation}
\mathcal{\vec{X}}_k = \mathcal{W}_k*\mathcal{D}^{-1}_k*\mathcal{W}^T_k * \mathcal{\vec{B}}.
\label{eq: 6}
\end{equation}
The approach so described for computing the approximate solution of \eqref{I3} is referred to as the truncated tEVD (T-tEVD) method. Analogously, the approximate solution of \eqref{I6} by the T-tEVD method is determined by
\begin{equation}
\mathcal{X}_k = \mathcal{W}_k*\mathcal{D}^{-1}_k*\mathcal{W}^T_k * \mathcal{B}.
\label{eq: 61}
\end{equation}
Our implementations of the solution methods for \eqref{eq: 6} and \eqref{eq: 61} first compute the full tEVD of $\mathcal{A}$, and then uses the discrepancy principle to decide where to truncate. The resulting approach based on $k$-term truncation is referred to as the T-tEVD method.

\subsection{The solution method for \eqref{I3} and \eqref{I6} by the t-Lanczos process}\label{sec3.2}

A popular approach for computing a few extreme eigenvalues and associated eigenvectors of a large symmetric matrix $A\in\mathbb{R}^{m \times m}$ are based on the symmetric Lanczos iteration; see, e.g., \cite{AGRR,GORR}. Analogously for third order tensors, the t-Lanczos process described by Algorithm \ref{t-Lanc} can be used to compute a few extreme eigentubes and associated eigenmatrices of a large symmetric tensor $\mathcal{A}\in \mathbb{R}^{m \times m \times n}$; see Kilmer et al. \cite{KBHH}.

The t-Lanczos process uses the normalization Algorithm \ref{normalize} described in \cite{KBH}. 

\vspace{.3cm}
\begin{algorithm}[H]
\SetAlgoLined
\KwIn{$\mathcal{\vec{X}} \in \mathbb{R}^{m \times 1 \times n} \neq \mathcal{\vec{O}}$}
\KwOut{$ \mathcal{\vec{V}}*\mathbf{a}$ with $\|\mathcal{\vec{V}}\| = 1$}
$\mathcal{\vec{V}} = \mathtt{fft}(\mathcal{\vec{X}},[\;],3)$\\
\For{$j = 1$ \bf{to} $n$}{
$\mathbf{a}^{(j)} \gets \|\mathcal{\vec{V}}^{(j)}\|_2 \;\;\;$ ($\mathcal{\vec{V}}^{(j)}$ is a vector)\\
\eIf{$\mathbf{a}^{(j)} > \mathtt{tol}$}{
$\mathcal{\vec{V}}^{(j)} \gets \frac{1}{\mathbf{a}^{(j)}} \mathcal{\vec{V}}^{(j)}$\\
}{
$\mathcal{\vec{V}}^{(j)} \gets \mathtt{randn}(m,1); \;\; \mathbf{a}^{(j)} \gets \|\mathcal{\vec{V}}^{(j)}\|_2;\;\; \mathcal{\vec{V}}^{(j)} \gets \frac{1}{\mathbf{a}^{(j)}} \mathcal{\vec{V}}^{(j)}; \;\; \mathbf{a}^{(j)} \gets 0$\\
}
}
$\mathcal{\vec{V}} \gets \mathtt{ifft}(\mathcal{\vec{V}},[\;],3); \;\; \mathbf{a} \gets \mathtt{ifft}(\mathbf{a},[\;],3)$
\caption{Normalize \cite{KBH}}
\label{normalize}
\end{algorithm}

Algorithm \ref{normalize} takes a nonzero tensor $\mathcal{\vec{X}} \in \mathbb{R}^{m \times 1 \times n}$, and returns $\mathcal{\vec{V}} \in \mathbb{R}^{m \times 1 \times n}$ and $\mathbf{a} \in \mathbb{R}^{1 \times 1 \times n}$
such that
\begin{equation*}
\mathcal{\vec{X}} = \mathcal{\vec{V}}* \mathbf{a} \;\;\; \text{and} \;\;\; \|\mathcal{\vec{V}}\| = 1,
\end{equation*}
where \[
\|\mathcal{\vec{X}}\| := \frac{\|\vec{\mathcal{X}}^T*\vec{\mathcal{X}}\|_F}{\|\vec{\mathcal{X}}\|_F},
\]
and $\|\mathcal{\vec{X}}\| = 0$ if $\mathcal{\vec{X}} = 0$. The tubal scalar $\mathbf{a}$ might not be invertible; see \cite{KBHH}. We remark that $\mathbf{a}$ is invertible if there exists $\mathbf{b}$ such that $\mathbf{a*b} = \mathbf{b*a} = {\bf e}_1$. In Algorithm \ref{normalize}, $\mathbf{a}^{(j)}$ is a scalar, whereas $\mathcal{\vec{V}}^{(j)}$ is a vector of length $m$. We remark that the Frobenius norm of a tensor column $\mathcal{\vec{X}}\in\mathbb{R}^{m\times 1\times n}$
given by
\begin{equation}\label{tfrobX}
\|\mathcal{\vec{X}}\|_F = \sqrt{\big(\mathcal{\vec{X}}^T*\mathcal{\vec{X}}\big)^{(1)}},
\end{equation}
is the square root of  the first frontal slice of the tube $\mathcal{\vec{X}}^T*\mathcal{\vec{X}}\in\mathbb{R}^{1\times 1\times n}$; see \cite{KBHH}.

Assume for simplicity that Algorithm \ref{t-Lanc} does not break down, i.e., choose $k$ small enough so that $\mathbf{c}_i$ is always invertible ($\mathbf{\widehat{c}}_i$ of $\mathbf{c}_i$ has no zero Fourier coefficients) for $1\leq i \leq k$. Then $k$ steps of Algorithm \ref{t-Lanc} with initial tensor $\mathcal{\vec{B}}$ will produce the t-Lanczos decomposition

\begin{equation}
\mathcal{A}*\mathcal{{Q}}_k = \mathcal{{Q}}_{k+1}*\mathcal{\bar{T}}_k,
\label{t-lanczos}
\end{equation}
where
\[
\mathcal{Q}_j = [\mathcal{\vec{Q}}_1, \mathcal{\vec{Q}}_2, \dots, \mathcal{\vec{Q}}_j] \in \mathbb{R}^{m \times j \times n}, ~~~ j \in \{k,k+1\},
\]
and the lateral slices $\mathcal{\vec{Q}}_j$, $j=1,2,\dots,k$, form an orthonormal tensor basis for the t-Krylov subspace \eqref{krylov2}. The initial tensor column $\mathcal{\vec{Q}}_1$ is obtained by normalizing $\mathcal{\vec{B}}$ such $\mathcal{\vec{B}} = \mathcal{\vec{Q}}_1*\mathbf{z}_0$ and $\|\mathcal{\vec{Q}}_1\| = 1$. Moreover,
\be \label{fra}
\mathcal{Q}_{k+1}^T * \mathcal{\vec{B}} = \mathcal{Q}_{k+1}^T*\mathcal{\vec{Q}}_1*\mathbf{z}_0 = \vec{e}_1 * \mathbf{z}_0,
\ee
and the tridiagonal tensor $\mathcal{T}_k$ is defined by
\begin{center}
$\mathcal{\bar{T}}_k = \begin{bmatrix}
\mathbf{c}_1 & \mathbf{z}_1& &\\
\mathbf{z}_1 & \mathbf{c}_2 & \mathbf{z}_2 \\
& \mathbf{z}_2 & \mathbf{c}_3 & \ddots\\
& & \ddots & \ddots & \mathbf{z}_{k-1} \\
& & & \mathbf{z}_{k-1} & \mathbf{c}_k\\
&&&&\mathbf{z}_k
\end{bmatrix} \in \mathbb{R}^{ (k+1)\times k\times n}$.
\end{center}

\vspace{.3cm}
\begin{algorithm}[H]
\SetAlgoLined
\KwIn{$\mathcal{A} \in \mathbb{R}^{m \times m \times n}$, $\mathcal{\vec{B}} \in \mathbb{R}^{m \times 1 \times n}$}
$\mathcal{\vec{Q}}_0 \leftarrow \mathcal{\vec{O}}$\\
$[\mathcal{\vec{Q}}_1, {\mathbf{z}}_0] \leftarrow \text{Normalize}(\mathcal{\vec{B}})$\\
\For {$i = 1$ \bf{to} $k$}{
$\mathcal{\vec{V}} \leftarrow \mathcal{A}*\mathcal{\vec{Q}}_i$\\
$\mathbf{c}_i \leftarrow \mathcal{\vec{Q}}^T_i*\mathcal{\vec{V}}$\\
$\mathcal{\vec{V}} \leftarrow \mathcal{\vec{V}} - \mathcal{\vec{Q}}_{i-1}*\mathbf{z}_{i-1} - \mathcal{\vec{Q}}_i * \mathbf{c}_i$\\
$[\mathcal{\vec{Q}}_{i+1}, {\mathbf{z}}_i] \leftarrow \text{Normalize}(\mathcal{\vec{V}})$
}
\caption{Symmetric t-product Lanczos (t-Lanczos) process \cite{KBHH}}
\label{t-Lanc}
\end{algorithm}\vspace{.3cm}

We use the t-Lanczos decomposition \eqref{t-lanczos} to reduce \eqref{I3} to a problem of small size. Let $\mathcal{\vec{X}} = \mathcal{Q}_k*\mathcal{\vec{Y}}$. Then by using \eqref{fra} and \eqref{normF}, the minimization problem \eqref{I3} reduces to

\begin{equation}
\min_{\mathcal{\vec{Y}} \in \mathbb{R}^{k\times 1 \times n}} \|\mathcal{\bar{T}}_k* \mathcal{\vec{Y}} - \mathcal{Q}_{k+1}^T*\mathcal{\vec{B}}\|_F = \min_{\mathcal{\vec{Y}} \in \mathbb{R}^{k \times 1 \times n}} \|\mathcal{\bar{T}}_k* \mathcal{\vec{Y}} - \vec{e}_1 * \mathbf{z}_0\|_F.
\label{eq: 9t}
\end{equation}
Let $\mathcal{\vec{Y}}_k$ be the solution of \eqref{eq: 9t} determined by using Algorithm \ref{backslash}. Then the approximate solution of \eqref{I3} is expressed as
\[
\mathcal{\vec{X}}_k = \mathcal{Q}_k * \mathcal{\vec{Y}}_k.
\]
We refer to this approach of solving \eqref{I3} as the t-Lanczos method. This method is implemented in the computed examples in Section \ref{sec5} with reorthogonalization of the t-Lanczos process.

\vspace{.3cm}
\begin{algorithm}[H]
\SetAlgoLined
\KwIn{$\mathcal{C}\in\mathbb{R}^{\ell\times m\times n}$, where its Fourier transform has
nonsingular frontal slices; $\mathcal{\vec{D}}\in\mathbb{R}^{\ell\times 1\times n}$,
$\mathcal{\vec{D}}\neq \mathcal{\vec{O}}$}
\KwOut{The solution $\mathcal{\vec{Y}}\in\mathbb{R}^{m\times 1\times n}$ of
$\min_{\mathcal{\vec{Y}}\in\mathbb{R}^{m\times 1\times n}}
\|\mathcal{C*\vec{Y}}-\mathcal{\vec{D}}\|_F$}
$\mathcal{{C}} \leftarrow \mathtt{fft}(\mathcal{C},[\;],3)$\\
$\mathcal{{\vec{D}}} \leftarrow \mathtt{fft}(\mathcal{\vec{D}},[\;],3)$\\
\For {$i =1$ \bf{to} $ n$ }{
$\mathcal{{\vec{Y}}}(:,:,i)=\mathcal{{C}}(:,:,i)\backslash\mathcal{{\vec{D}}}(:,:,i)$,
where $\backslash$ denotes MATLAB's backslash operator
}
$\mathcal{\vec{Y}} \leftarrow \mathtt{ifft}(\mathcal{{\vec{Y}}},[\;],3)$
\caption{Solution of a generic tensor least squares problem \cite{RU1}}
\label{backslash}
\end{algorithm}
\vspace{.3cm}

The approximate solution of \eqref{I6} is determined by applying the t-Lanczos method to solve each one of the $p$ minimization problems
\begin{equation}\label{33no}
\min_{\mathcal{\vec{X}}_j \in \mathbb{R}^{m\times 1\times n}}
\|\mathcal{A}*\mathcal{\vec{X}}_j-\mathcal{\vec{B}}_j\|_F, ~~~ j = 1,2,\dots,p, ~~~ p>1,
\end{equation}
where $\mathcal{\vec{B}}_1,\mathcal{\vec{B}}_2,\dots,\mathcal{\vec{B}}_p$ are the lateral slices of the data tensor $\mathcal{B}\in \mathbb{R}^{m\times p \times n}$. For each problem \eqref{33no}, we let $\mathcal{\vec{B}}_{j,{\rm true}}$ denote the unknown error-free tensor slice associated
with the available error-contaminated tensor slice $\mathcal{\vec{B}}_j$, and assume that bounds
\[
\|\mathcal{\vec{E}}_j\|_F\leq\delta_j,\quad j=1,2,\ldots,p,
\]
for the Frobenius norm of the errors,
\[
\mathcal{E}_j:=\mathcal{\vec{B}}_j-\mathcal{\vec{B}}_{j,{\rm true}},\quad j=1,2,\dots,p,
\]
are available or can be estimated; cf. \eqref{errbd}. Then the solution method for \eqref{I6} obtained by applying the t-Lanczos method to each one of the $p$ problems \eqref{33no} is referred to as the t-Lanczos$_p$ method.

We conclude this subsection by describing how an approximate solution of \eqref{I3} can be expressed in terms of the t-eigenpairs of largest magnitude. Using the t-Lanczos decomposition \eqref{t-lanczos}, we have
\be \label{slanc}
\mathcal{{Q}}_k^T*\mathcal{A}*\mathcal{{Q}}_k = \mathcal{{T}}_k,
\ee
where $\mathcal{{T}}_k$ is the $k \times k \times n$ leading sub-tridiagonal tensor of $\mathcal{\bar{T}}_k$. Moreover, since $k$ is small, we can compute the full tEVD of $\mathcal{T}_k$,
\begin{equation}
\mathcal{T}_k = \mathcal{W} * \mathcal{D} * \mathcal{W}^T,
\label{eq: 8}
\end{equation}
where $\mathcal{W} = [\mathcal{\vec{W}}_1, \mathcal{\vec{W}}_2, \dots, \mathcal{\vec{W}}_k] \in \mathbb{R}^{k\times k\times n}$ has orthonormal tensor columns, and
\[
\mathcal{D} = \text{diag}[{\mathbf{d}}_1, {\mathbf{d}}_2, \dots, {\mathbf{d}}_k] \in \mathbb{R}^{k\times k\times n}
\]
is an f-diagonal tensor with eigentubes ${\mathbf{d}}_j$ corresponding to the eigenmatrices $\mathcal{\vec{W}}_j$, $j=1,2,\ldots, k$.

Substituting $\mathcal{\vec{X}} = \mathcal{Q}_k*\mathcal{\vec{Y}}$ in \eqref{I3}, using \eqref{slanc}, \eqref{eq: 8}, and \eqref{normF} with $\mathcal{\vec{Z}} = \mathcal{W}^T*\mathcal{\vec{Y}}$, gives
\begin{equation}
\min_{\mathcal{\vec{Z}}\in \mathbb{R}^{k\times 1\times n}} \|\mathcal{D}* \mathcal{\vec{Z}} - (\mathcal{Q}_k*\mathcal{W})^T*\mathcal{\vec{B}}\|_F.
\label{eq: 9}
\end{equation}
The solution of \eqref{eq: 9} is given by $\mathcal{\vec{Z}}_k = \mathcal{D}^{-1}*(\mathcal{Q}_k*\mathcal{W})^T * \mathcal{\vec{B}}$. Hence, an approximate solution of \eqref{I3} can be expressed as
\begin{equation}
\mathcal{\vec{X}}_k = (\mathcal{Q}_k * \mathcal{W})*\mathcal{D}^{-1}*(\mathcal{Q}_k*\mathcal{W})^T * \mathcal{\vec{B}}.
\label{eq: 10}
\end{equation}
The number of steps required by the t-Lanczos process is determined by the discrepancy principle \eqref{disp}. Similarly, by replacing $\mathcal{\vec{B}}$ in \eqref{eq: 10} with $\mathcal{B}$, we obtain an approximate solution $\mathcal{X}_k$ of \eqref{I6}.

Analogously to matrix computations, the lateral slices of $\mathcal{Q}_k * \mathcal{W}$ are called the Ritz matrices \cite{KBHH}. We refer to the pairs $(\mathbf{d}_i,\mathcal{Q}_k*\mathcal{\vec{W}}_i)$ for $1 \leq i \leq k $ as the t-eigenpairs of largest magnitude. They correspond to the eigentubes of largest Frobenius norm and associated eigenmatrices.

\section{Nonsymmetric t-linear discrete ill-posed problems}\label{sec4}
This section is concerned with the solution for \eqref{I3} and \eqref{I6} when the tensor $\mathcal{A}\in \mathbb{R}^{\ell \times m\times n}$ is not necessarily symmetric. We consider the solution of both problems by using the T-tSVD \cite{KM}, and the tGKB process \cite{KBHH}. A new approach for determining the approximate solution of \eqref{I3} and \eqref{I6} based on randomized T-tSVD (RT-tSVD) is presented in Section \ref{sec4.2}.

\subsection{The T-tSVD method for the approximate solution of \eqref{I3} and \eqref{I6}}\label{sec4.1}
This subsection describes the solution method for \eqref{I3} and \eqref{I6} by using the T-tSVD introduced in \cite{KM}. Let the tSVD of $\mathcal{A} \in \mathbb{R}^{\ell \times m\times n}$, $\ell \geq m$, be given by
\begin{equation}
\mathcal{A} = \mathcal{U}* \mathcal{S} *\mathcal{V}^T \implies \mathcal{A}*\mathcal{\vec{V}}_j = \mathcal{\vec{U}}_j*{\mathbf{s}}_j, ~~ j = 1, 2,\ldots, m,
\label{eq: 11}
\end{equation}
where $\mathcal{U} = [\mathcal{\vec{U}}_1, \mathcal{\vec{U}}_2, \dots, \mathcal{\vec{U}}_m] \in \mathbb{R}^{\ell \times m \times n}$ is partially orthogonal, $\mathcal{V} = [\mathcal{\vec{V}}_1, \mathcal{\vec{V}}_2, \dots, \mathcal{\vec{V}}_m] \in \mathbb{R}^{m \times m \times n}$ is orthogonal, and $\mathcal{S} = \text{diag}[{\mathbf{s}}_1, {\mathbf{s}}_2, \dots, {\mathbf{s}}_m] \in \mathbb{R}^{m \times m \times n}$ is an f-diagonal tensor with singular tubes ${\mathbf{s}}_j$, $j=1,2,\dots,m$, corresponding to the left and right singular matrices $\mathcal{\vec{U}}_j$ and $\mathcal{\vec{V}}_j$, respectively. The Frobenius norm of the singular tubes ${\mathbf{s}}_j$ are assumed to decay rapidly, and are ordered according to
\[
\|{\mathbf{s}}_1\|_F \geq \|{\mathbf{s}}_2\|_F \geq \dots \geq \|{\mathbf{s}}_k\|_F > \|{\mathbf{s}}_{k+1}\|_F = \dots = \|{\mathbf{s}}_p\|_F = 0,
\]
where $k$ is the tubal rank of $\mathcal{A}$. Let
\[
\mathcal{U}_s = [\mathcal{\vec{U}}_1, \mathcal{\vec{U}}_2, \dots, \mathcal{\vec{U}}_s] \in \mathbb{R}^{\ell \times s\times n}, ~\mathcal{V}_s = [\mathcal{\vec{V}}_1, \mathcal{\vec{V}}_2, \dots, \mathcal{\vec{V}}_s] \in \mathbb{R}^{m\times s\times n},~\mathcal{S}_s = {\rm diag}[{\mathbf{s}}_1, {\mathbf{s}}_2, \dots, {\mathbf{s}}_s] \in \mathbb{R}^{s\times s\times n}.
\]

Define the T-tSVD of $\mathcal{A}$ as
\begin{equation}
\mathcal{A}_s = \mathcal{U}_s* \mathcal{S}_s *\mathcal{V}^T_s = \sum_{i=1}^s \mathcal{\vec{U}}_i*\mathbf{s}_i*\mathcal{\vec{V}}_i^T,
\label{eq: 12}
\end{equation}
for some $1 \leq s \leq k$. Then by Theorem \ref{thm 2.1}, $\mathcal{A}_s$ is the best tubal rank-$s$ approximation to $\mathcal{A}$. The approximate solution of \eqref{I3} by the T-tSVD is less contaminated by propagated error since the condition numbers of $\mathcal{A}_s^{(i)}$, $i = 1,2,\dots,n$, are smaller than those of $\mathcal{A}^{(i)}$.

Replacing $\mathcal{A}$ by $\mathcal{A}_s$ and $\mathcal{\vec{B}}$ by $\mathcal{U}_s* \mathcal{U}_s^T * \mathcal{\vec{B}}$ (the orthogonal projector onto the range of $\mathcal{A}_s$) in \eqref{I3} for a suitably small $s$, setting $\mathcal{\vec{X}} = \mathcal{V}_s * \mathcal{\vec{Y}}$ for some $ \mathcal{\vec{Y}}\in \mathbb{R}^{s\times 1\times n}$, and using \eqref{normF} gives the minimization problem
\begin{equation}
\min_{\mathcal{\vec{Y}}\in \mathbb{R}^{s\times 1\times n}} ||\mathcal{S}_s * \mathcal{\vec{Y}} - \mathcal{U}^T_s * \mathcal{\vec{B}}||_F.
\label{eq: 13}
\end{equation}
The solution of \eqref{eq: 13} is given by $\mathcal{\vec{Y}}_s = \mathcal{S}^{-1}_s*\mathcal{U}^T_s * \mathcal{\vec{B}}$. Hence, an approximate solution of \eqref{I3} is expressed as
\begin{equation}
\mathcal{\vec{X}}_s = \mathcal{V}_s*\mathcal{S}^{-1}_s*\mathcal{U}^T_s * \mathcal{\vec{B}} = \sum_{i=1}^s \mathcal{\vec{V}}_i*\mathbf{s}_i^{-1}*\mathcal{\vec{U}}_i^T * \mathcal{\vec{B}}, \;\; \;\;\;\;1 \leq s \leq k.
\label{eq: 14}
\end{equation}
We refer to this approach of computing an approximate solution of \eqref{I3} as the T-tSVD method. The implementation of the T-tSVD method proceeds by first computing a full tSVD of $\mathcal{A}$ and then determining the truncation index $s$ by the discrepancy principle \eqref{disp}. The T-tSVD method has the advantage that the computation of the T-tSVD of $\mathcal{A}$ is independent of the data slice $\mathcal{\vec{B}}$; see, e.g., \cite{NR} for discussions on the truncated SVD for ill-posed problems in two dimensions. We remark that the approximate solution of \eqref{I6} follows analogously by replacing the data slice $\mathcal{\vec{B}}$ by a general third order tensor $\mathcal{B}\in \mathbb{R}^{\ell\times p\times n}$, $p>1$.

\subsection{The randomized T-tSVD method for the solution of \eqref{I3} and \eqref{I6}}\label{sec4.2}

This subsection describes the randomized T-tSVD (RT-tSVD) method for the approximate solution of linear discrete ill-posed problems \eqref{I3} and \eqref{I6}. The randomized tSVD (R-tSVD) method with the t-product was first introduced by Zhang et al. \cite{ZSKA} and applied to facial recognition problems. Here we focus on ill-posed problems arising in image and video restorations.

We present a new R-tSVD method that is based on the randomized block algorithms for efficient low-rank matrix computations described in \cite{MV, YGL}. The proposed R-tSVD method is described by Algorithm \ref{Rtsvd}. It determines small-size tensors $\mathcal{Q}\in \mathbb{R}^{\ell \times r \times n}$, with $r$ orthonormal tensor columns that form a basis for the range of $\mathcal{A}$, and $\mathcal{\tilde{B}}\in \mathbb{R}^{r \times m\times n}$ such that
\[
\mathcal{A} \approx \mathcal{Q}*\mathcal{\tilde{B}},
\]
where $\mathcal{\tilde{B}} = \mathcal{Q}^T*\mathcal{A}$, and
\be \label{apperr}
\|\mathcal{A} - \mathcal{Q}*\mathcal{\tilde{B}} \|_F^2 = \|\mathcal{A}\|_F^2 - \|\mathcal{\tilde{B}} \|_F^2 < \epsilon^2,
\ee
for some tolerance $\epsilon>0$. A full tSVD of $\mathcal{\tilde{B}}$ is then computed to obtain a low tubal rank approximation of $\mathcal{A}$,
\[
\mathcal{A} \approx \mathcal{Q}*\mathcal{\tilde{B}} = (\mathcal{Q}*\mathcal{\tilde{U}})*\mathcal{\tilde{S}}*\mathcal{\tilde{V}}^T,
\]
where $\mathcal{\tilde{B}} = \mathcal{\tilde{U}}*\mathcal{\tilde{S}}*\mathcal{\tilde{V}}^T$, $\mathcal{\tilde{U}}\in \mathbb{R}^{r \times r \times n}$, $\mathcal{\tilde{S}}\in \mathbb{R}^{r \times m \times n}$, and $\mathcal{\tilde{V}}\in \mathbb{R}^{m \times m \times n}$. The equality in \eqref{apperr} is shown in Theorem \ref{the1}. 

The computation of an accurate approximate factorization of $\mathcal{A}$, of tubal rank-$k$, that is close to the optimal tubal rank by Algorithm \ref{Rtsvd}, can be achieved by using a Gaussian random tensor slice $\mathcal{\vec{G}}\in \mathbb{R}^{m \times 1 \times n}$ at each $r$ step such that
\[ r = k + \rho, \]
where $\rho$ is a small oversampling parameter, say, $\rho=5$ or $\rho=10$, and $k$ is determined by the discrepancy principle.

Denote the low tubal rank approximation of $\mathcal{A}$ by $\mathcal{A}_k$. Let
\[\mathcal{U}_k = \mathcal{Q}*\mathcal{\tilde{U}}(:,1:k,:), ~\mathcal{S}_k = \mathcal{\tilde{S}}(1:k,1:k,:), {\rm and}~ \mathcal{V}_k = \mathcal{\tilde{V}}(:,1:k,:).\] Then
\[
\mathcal{A}_k= \mathcal{U}_k * \mathcal{S}_k* \mathcal{V}_k^T,
\]
where $\mathcal{U}_k$ and $\mathcal{V}_k$ have $k$ orthonormal tensor columns, and $\mathcal{S}_k$ is f-diagonal containing $k$ singular tubes of the largest Frobenius norm. We use the representation above to determine an approximate solution of \eqref{I3} and \eqref{I6} by following a similar approach described in Section \ref{sec4.1}. Hence,
\[
\mathcal{\vec{X}}_k = \mathcal{V}_k*\mathcal{S}_k^{-1}*\mathcal{U}_k^T*\mathcal{\vec{B}} ~~~ {\rm and}~~~ \mathcal{X}_k = \mathcal{V}_k*\mathcal{S}_k^{-1}*\mathcal{U}_k^T*\mathcal{B},
\]
respectively. We refer to the described method above as the randomized T-tSVD (RT-tSVD) method. Similar to the T-tSVD and T-tEVD methods, a low tubal rank factorization of $\mathcal{A}$ determined by the RT-tSVD method is independent of the data tensors $\mathcal{\vec{B}}$ and $\mathcal{B}$. 

\vspace{.3cm}
\begin{algorithm}[H]
\SetAlgoLined
\KwIn{$ \mathcal{A} \in \mathbb{R}^{\ell \times m \times n}$, $\mathcal{\vec{Q}}_0 = \mathcal{\vec{O}}\in \mathbb{R}^{\ell \times 1 \times n}$, $\mathcal{\tilde{B}}(0,:,:) = \mathcal{\vec{O}}\in \mathbb{R}^{1\times m \times n}$, $r=0$, $\epsilon$}
\KwOut{$\mathcal{\tilde{U}} \in \mathbb{R}^{\ell \times r \times n}$, $\mathcal{\tilde{S}} \in \mathbb{R}^{r \times m \times n}$, $\mathcal{\tilde{V}} \in \mathbb{R}^{m \times m \times n}$}
$\eta \leftarrow \|\mathcal{A} \|_F^2$\\
\While{$\eta \geq \epsilon^2$}{
$r \leftarrow r+1$\\
Generate a Gaussian random tensor $\mathcal{\vec{G}} \in \mathbb{R}^{m \times 1 \times n}$\\
$\mathcal{\vec{Z}} \leftarrow \mathcal{A*\vec{G}} - \mathcal{\vec{Q}}_{r-1}*\mathcal{\tilde{B}}(r-1,:,:)*\mathcal{\vec{G}}$ \\
$[\mathcal{\vec{Z}},\sim] \leftarrow {\tt Normalize}(\mathcal{\vec{Z}})$ by Algorithm \ref{normalize}\\
$\mathcal{\vec{Z}} \leftarrow \mathcal{\vec{Z}} - \sum\limits_{j=1}^{r-1} \mathcal{\vec{Q}}_j*(\mathcal{\vec{Q}}_j^T * \mathcal{\vec{Z}})$\\
$[\mathcal{\vec{Z}},\sim] \leftarrow {\tt Normalize}(\mathcal{\vec{Z}})$ by Algorithm \ref{normalize}\\
$\mathcal{\vec{Q}}_r \leftarrow \mathcal{\vec{Z}}$\\
$\mathcal{\vec{Y}}\leftarrow \mathcal{\vec{Z}}^T*\mathcal{A}$\\
$\mathcal{\tilde{B}}(r,:,:) \leftarrow \mathcal{\vec{Y}}$\\
$\eta \leftarrow \eta - \|\mathcal{\vec{Y}} \|_F^2$\\
}
Form $\mathcal{Q} \leftarrow [\mathcal{\vec{Q}}_1, \mathcal{\vec{Q}}_2, \dots, \mathcal{\vec{Q}}_r] \in \mathbb{R}^{\ell \times r \times n}$ and $\mathcal{\tilde{B}} \leftarrow \begin{bmatrix}\mathcal{\tilde{B}}(1,:,:) \\ \vdots \\ \mathcal{\tilde{B}}(r,:,:) \end{bmatrix} \in \mathbb{R}^{r \times m \times n}$\\
Compute the full tSVD of $\mathcal{\tilde{B}}$ to obtain $\mathcal{\breve{U}}\in \mathbb{R}^{r \times r \times n}$, $\mathcal{\tilde{S}}\in \mathbb{R}^{r \times m \times n}$ and $\mathcal{\tilde{V}}\in \mathbb{R}^{m \times m \times n}$\\
Compute $\mathcal{\tilde{U}} \leftarrow \mathcal{Q*\breve{U}}$\\
\caption{Randomized tSVD (R-tSVD)}
\label{Rtsvd}
\end{algorithm}\vspace{.3cm}

Algorithm \ref{Rtsvd} differs from the one presented in \cite{ZSKA} because it iteratively builds the lateral and horizontal slices of $\mathcal{Q}$ and $\mathcal{\tilde{B}}$, respectively, and measures the approximation error by using the right-hand side of \eqref{apperr}. The justifications for the latter choice of the error indicator in step $12$ of Algorithm \ref{Rtsvd} are presented in Theorems $4.1$ and $4.2$; see \cite{YGL} for discussions on the matrix case. The advantage of this approach is that the Frobenius norm of the residual can be computed once $\mathcal{\tilde{B}}$ is accessible and $\|\mathcal{A} \|_F^2$ is known {\it a prior}. In step $7$ of Algorithm \ref{Rtsvd}, reorthogonalization is carried out to ensure that the lateral slices $\mathcal{\vec{Q}}_j$, $j=1,2,\dots,r$, of $\mathcal{Q}$ are orthonormal.

\begin{thm}\label{the1}
Let $\mathcal{A}\in \mathbb{R}^{\ell \times m \times n}$, and let $\mathcal{Q}\in \mathbb{R}^{\ell \times r \times n}$ be a partially orthogonal tensor such that $\mathcal{\tilde{B}} = \mathcal{Q}^T*\mathcal{A}$. Then
\[
\|\mathcal{A} - \mathcal{Q}*\mathcal{\tilde{B}} \|_F^2 = \|\mathcal{A}\|_F^2 - \|\mathcal{\tilde{B}} \|_F^2.
\]
Proof: Recall that
\[
\|\mathcal{A}\|_F^2 = {\tt trace}\left( \left( \mathcal{A}^T*\mathcal{A} \right)^{(1)} \right).
\]
Then
\begin{equation*}
\begin{split}
(\mathcal{A} - \mathcal{Q}*\mathcal{\tilde{B}})^T(\mathcal{A} - \mathcal{Q}*\mathcal{\tilde{B}}) &= (\mathcal{A} - \mathcal{Q}*\mathcal{Q}^T*\mathcal{A})^T(\mathcal{A} - \mathcal{Q}*\mathcal{Q}^T*\mathcal{A}) \\ &= \mathcal{A}^T *\mathcal{A} - 2\mathcal{A}^T*\mathcal{Q}*\mathcal{Q}^T*\mathcal{A} + \mathcal{A}^T*\mathcal{Q}*\mathcal{Q}^T *\mathcal{Q}*\mathcal{Q}^T*\mathcal{A}\\ &= \mathcal{A}^T *\mathcal{A} - \mathcal{A}^T*\mathcal{Q}*\mathcal{Q}^T*\mathcal{A} \\ &= \mathcal{A}^T *\mathcal{A} - \mathcal{\tilde{B}}^T *\mathcal{\tilde{B}}.
\end{split}
\end{equation*}
Hence,
\be \label{trace}
\left((\mathcal{A} - \mathcal{Q}*\mathcal{\tilde{B}})^T(\mathcal{A} - \mathcal{Q}*\mathcal{\tilde{B}})\right)^{(1)} = \left(\mathcal{A}^T *\mathcal{A}\right)^{(1)} - \left(\mathcal{\tilde{B}}^T *\mathcal{\tilde{B}}\right)^{(1)}.
\ee
The proof follows by applying the {\tt trace} operator to both sides of \eqref{trace}. ~~~$\Box$
\end{thm}
The following theorem shows that the error indicator $\eta$ is the square of the approximation error, and it is updated by computing the Frobenius norm of $\mathcal{\tilde{B}}(r,:,:)$.
\begin{thm}\label{the2}
Let $\eta^{[r]}$, $\mathcal{Q}^{[r]}$, and $\mathcal{\tilde{B}}^{[r]}$ correspond to the terms of $\eta$, $\mathcal{Q}$, and $\mathcal{\tilde{B}}$ determined by Algorithm \ref{Rtsvd} at the $r$-th iteration, respectively. Then
\[
\eta^{[r]} = \|\mathcal{A} -\mathcal{Q}^{[r]}*\mathcal{\tilde{B}}^{[r]} \|^2_F.
\]
Proof: Using step $12$ of Algorithm \ref{Rtsvd} and the property of Frobenius norm, we have
\[
\eta^{[r]} = \|\mathcal{A}\|^2_F - \sum_{i=1}^r \| \mathcal{B}(i,:,:)\|^2_F = \|\mathcal{A}\|_F^2 - \|\mathcal{\tilde{B}}^{[r]} \|_F^2.
\]
Since $\mathcal{Q}^{[r]}$ has orthonormal tensor columns, and $\mathcal{\tilde{B}}^{[r]} = \left(\mathcal{Q}^{[r]}\right)^T*\mathcal{A}$, the proof follows by Theorem \ref{the1}. ~~~ $\Box$
\end{thm}

We conclude this section with the error bound for the R-tSVD method. It can be shown that the expected error of the probabilistic part of the R-tSVD method is given by
\begin{equation*}
\mathbb{E}\|\mathcal{A} - \mathcal{Q}*\mathcal{Q}^T*\mathcal{A} \|_F^2 \leq \Bigg(1 + \frac{k}{(k-\rho)}\Bigg) \Bigg(\frac{1}{n} \sum_{i=1}^n \sum_{j=k+1}^{\min\{\ell,m\}} \widehat{\sigma}_j^{(i)} \Bigg),
\end{equation*}
for $\rho \geq 2$ where $\mathbb{E}$ denotes expectation, and $\widehat{\sigma}_j^{(i)}$ is the $i$th component of $\mathcal{\widehat{S}}(j,j,:)$; see \cite{ZSKA}. We remark that the expected error for the R-tSVD method is within a factor of $1 + \frac{k}{(k-\rho)}$ which is the same as in the matrix case described by Halko et al. \cite{HMT}. The performance of the R-tSVD method also depends on the rate of decay of the singular values  of $\mathcal{\widehat{A}}^{(i)}$.

\subsection{The solution method for \eqref{I3} and \eqref{I6} by the tGKB process}\label{sec4.3}
This subsection reviews the t-product Golub-Kahan bidiagonalization (tGKB) process introduced by Kilmer et al. \cite{KBHH} for the computation of a few extreme singular tubes and associated left and right singular matrices of a large nonsymmetric tensor $\mathcal{A}\in \mathbb{R}^{\ell \times m \times n}$ with an initial tensor slice $\mathcal{\vec{B}}\in\mathbb{R}^{\ell\times 1\times n}$. The tGKB process is described by Algorithm \ref{tgkb} with reorthogonalization applied in steps $4$ and $6$.

The solution methods for \eqref{I3} and \eqref{I6} considered in \cite{RU1} combines the tGKB process with Tikhonov regularization. Here, we consider regularization of \eqref{I3} and \eqref{I6} by truncated tGKB iterations. We assume that the number of iterations of the tGKB process is small enough to avoid break down; see \cite{RU1}. Then Algorithm \ref{tgkb} after $k$ steps produces the tGKB decomposition
\begin{equation}\label{t-gkb}
\mathcal{A} * \mathcal{W}_k = \mathcal{Q}_{k+1} * \mathcal{\bar{P}}_k,
\end{equation}
where
\[\mathcal{Q}_{k+1} := [\mathcal{\vec{Q}}_1, \mathcal{\vec{Q}}_2, \dots, \mathcal{\vec{Q}}_{k+1}]\in \mathbb{R}^{ \ell \times {(k+1)}\times n}, ~~~\mathcal{W}_k := [\mathcal{\vec{W}}_1, \mathcal{\vec{W}}_2, \dots, \mathcal{\vec{W}}_k]\in \mathbb{R}^{ m\times k\times n},\]
and
\begin{center}
$\mathcal{\bar{P}}_k = \begin{bmatrix}
\mathbf{c}_1 & & &\\
\mathbf{z}_2 & \mathbf{c}_2 & \\
& \mathbf{z}_3 & \mathbf{c}_3 & \\
& & \ddots & \ddots & \\
& & & \mathbf{z}_{k} & \mathbf{c}_k\\
& & & & \mathbf{z}_{k+1}
\end{bmatrix} \in \mathbb{R}^{ (k+1)\times k\times n}$
\end{center}
is a lower bidiagonal tensor. The tensor columns $\mathcal{\vec{Q}}_i \in \mathbb{R}^{ \ell \times 1 \times n}$, $i=1,2,\dots, k+1$, and $\mathcal{\vec{W}}_i \in \mathbb{R}^{ m \times 1 \times n}$, $i=1,2,\dots, k$, generated by Algorithm \ref{tgkb} are orthonormal tensor bases for the t-Krylov subspaces $\mathbb{K}_{k+1}(\mathcal{A}*\mathcal{A}^T, \mathcal{\vec{B}})$ and $\mathbb{K}_k(\mathcal{A}^T*\mathcal{A}, \mathcal{A}^T*\mathcal{\vec{B}})$, respectively.

\vspace{.3cm}
\begin{algorithm}[H]
\SetAlgoLined
\KwIn{$\mathcal{A}\in\mathbb{R}^{\ell\times m\times n},\;
\mathcal{\vec{B}}\in\mathbb{R}^{\ell\times 1\times n}$}
$\mathcal{\vec{W}}_0 \leftarrow \mathcal{\vec{O}}$\\
$[\mathcal{\vec{Q}}_1, {\mathbf{z}}_1] \leftarrow \mathtt{Normalize}(\mathcal{\vec{B}})$
with $\mathbf{z}_1$ invertible\\
\For {$i = 1,2, \dots, k$}{
$\mathcal{\vec{W}}_i \leftarrow \mathcal{A}^T*\mathcal{\vec{Q}}_i - \mathcal{\vec{W}}_{i-1}* \mathbf{z}_i$\\
$\mathcal{\vec{W}}_i \leftarrow \mathcal{\vec{W}}_i - \sum\limits_{j=1}^{i-1} \mathcal{\vec{W}}_j*(\mathcal{\vec{W}}_j^T * \mathcal{\vec{W}}_i) \; (\text{reorthogonalization step})$\\
$[\mathcal{\vec{W}}_i, \mathbf{c}_i] \leftarrow \mathtt{Normalize}(\mathcal{\vec{W}}_i)$\\
$\mathcal{\vec{Q}}_{i+1} \leftarrow \mathcal{A}*\mathcal{\vec{W}}_i - \mathcal{\vec{Q}}_i* \mathbf{c}_i$ \\
$\mathcal{\vec{Q}}_{i+1} \leftarrow \mathcal{\vec{Q}}_{i+1} - \sum\limits_{j=1}^{i} \mathcal{\vec{Q}}_j*(\mathcal{\vec{Q}}_j^T * \mathcal{\vec{Q}}_{i+1})\; (\text{reorthogonalization step})$\\

$[\mathcal{\vec{Q}}_{i+1}, \mathbf{z}_{i+1}] \leftarrow \mathtt{Normalize}(\mathcal{\vec{Q}}_{i+1})$\\
}
\caption{Partial t-product Golub-Kahan bidiagonalization (tGKB) process \cite{RU1}}\label{tgkb}
\end{algorithm}\vspace{.3cm}

In the sequel, we describe how the tGKB decomposition \eqref{t-gkb} is used to determine an approximate solution of \eqref{I3}. Then we show that this solution in terms of the largest t-singular triplets of $\mathcal{A}$. The solution of \eqref{I6} by the tGKB process follows analogously by solving separately $p$ minimization problems \eqref{33no}.

Let $\mathcal{\vec{X}} = \mathcal{W}_k*\mathcal{\vec{Y}}$ for some tensor slice $\mathcal{\vec{Y}} \in \mathbb{R}^{k\times 1\times n}$. Using \eqref{t-gkb} and \eqref{normF}, the minimization problem \eqref{I3} reduces to
\begin{equation}\label{9tt}
\min_{\mathcal{\vec{Y}} \in \mathbb{R}^{k\times 1\times n}} \|\mathcal{\bar{P}}_k* \mathcal{\vec{Y}} - \mathcal{Q}_{k+1}^T*\mathcal{\vec{B}}\|_F = \min_{\mathcal{\vec{Y}} \in \mathbb{R}^{k\times 1\times n}} \|\mathcal{\bar{P}}_k* \mathcal{\vec{Y}} - \vec{e}_1 * \mathbf{z}_1\|_F,
\end{equation}
where $\mathcal{Q}_{k+1}^T*\mathcal{\vec{B}} = \vec{e}_1 * \mathbf{z}_1$ follows from Algorithm \ref{normalize}. Denote the solution of \eqref{9tt} by $\mathcal{\vec{Y}}_k$. Then the solution of \eqref{I3} is given by
\be\label{solgkb}
\mathcal{\vec{X}}_k = \mathcal{W}_k*\mathcal{\vec{Y}}_k.
\ee
We remark that the solution $\mathcal{\vec{Y}}_k$ of \eqref{9tt} is computed by using Algorithm \ref{backslash} until the discrepancy principle \eqref{disp} is satisfied. Reorthogonalization is not required by the tGKB method since only a few steps are needed.

Moreover, for sufficiently small $k$, we can compute the tSVD of $\mathcal{\bar{P}}_k$,
\begin{equation}
\mathcal{\bar{P}}_k = \mathcal{U} * \mathcal{S} * \mathcal{V}^T,
\label{eq: 16}
\end{equation}
where $\mathcal{U}\in \mathbb{R}^{(k+1)\times k\times n}$ has $k$ orthonormal tensor columns $\mathcal{\vec{U}}_j$, $j=1,2,\dots,k$, $\mathcal{V}\in \mathbb{R}^{k\times k\times n}$ is an orthogonal tensor with lateral slices $\mathcal{\vec{V}}_j$, and $\mathcal{S}\in \mathbb{R}^{k\times k\times n}$ is an f-diagonal tensor with singular tubes ${\mathbf{s}}_j$ corresponding to the left and right singular matrices $\mathcal{\vec{U}}_j$ and $\mathcal{\vec{V}}_j$, respectively. Substituting \eqref{eq: 16} into the left-hand side of \eqref{9tt}, and using \eqref{normF} gives the minimization problem
\begin{equation}
\min_{\mathcal{\vec{Y}}\in \mathbb{R}^{k\times 1\times n}} ||\mathcal{S}* \mathcal{\vec{Y}} - (\mathcal{Q}_{k+1} * \mathcal{U})^T*\mathcal{\vec{B}}||_F,
\label{eq: 17}
\end{equation}
with solution given by $\mathcal{\vec{Y}}_k = \mathcal{S}^{-1}*(\mathcal{Q}_{k+1} * \mathcal{U})^T * \mathcal{\vec{B}}$. Hence, the solution \eqref{solgkb} can equivalently be expressed in terms of the largest singular triplets as
\begin{equation}
\mathcal{\vec{X}}_k = (\mathcal{W}_k * \mathcal{V})*\mathcal{S}^{-1}*(\mathcal{Q}_{k+1}*\mathcal{U})^T * \mathcal{\vec{B}},
\label{solgkbb}
\end{equation}
where $\mathcal{S}^{-1}$ exist since we assume $k$ is chosen small to ensure the tGKB process does not break down.
The triplets $(s_i, \mathcal{Q}_{k+1}*\mathcal{\vec{U}}_i, \mathcal{W}_k*\mathcal{\vec{V}}_i)$ for $1 \leq i \leq k$ are analogues of the Ritz triplets in matrix computations \cite{KBHH}. They are the largest t-singular triplets of $\mathcal{A}$ corresponding to the singular tubes of largest Frobenius norm and associated left and right singular matrices, respectively. The solution $\mathcal{X}_k$ of \eqref{I6} in terms of the t-singular triplets replaces $\mathcal{\vec{B}}$ in \eqref{solgkbb} by $\mathcal{B}:=[\mathcal{\vec{B}}_1, \mathcal{\vec{B}}_2, \dots, \mathcal{\vec{B}}_p]\in \mathbb{R}^{\ell \times p\times n}$, $p>1$.

Finally, we describe an approach for the solution of \eqref{I6} based on t-Krylov recycling. We determine an approximate solution of \eqref{I6} by generating a t-Krylov subspace $\mathbb{K}_k(\mathcal{A}^T*\mathcal{A}, \mathcal{A}^T*\mathcal{\vec{B}}_1)$ of sufficiently large dimension $k$ to contain the solution of all the $p$ problems \eqref{33no}. Firstly, the approximate solution of the least squares problems \eqref{33no} for $p=1$ is determined in this t-Krylov subspace. Then using the same t-Krylov subspace, we determine an approximate solution of \eqref{33no} for $p=2$. However, the discrepancy principle may not be satisfied. Hence we increase the dimension $k$, and generate a new t-Krylov subspace $\mathbb{K}_k(\mathcal{A}^T*\mathcal{A}, \mathcal{A}^T*\mathcal{\vec{B}}_2)$ until the discrepancy principle is satisfied. After solving \eqref{33no} for $p=2$, we proceed analogously to solve \eqref{33no} for $p=3,4,\dots,p$. Reorthogonalization of the t-GKB process is required to ensure that
\[
\mathcal{Q}_{k+1}*\mathcal{\vec{B}}_j, ~~~j = 1,2,\dots,p,
\]
are computed with sufficient accuracy. This approach of solving \eqref{I6} thereby solving \eqref{33no} by nested t-Krylov subspaces is referred to as the {\tt nested}$\_$tGKB$_p$ method. This method is described by Algorithm \ref{nestKry} below. We remark that in situations when the bidiagonalization steps $k$ is large, restarting Algorithm \ref{nestKry} with a new tensor $\mathcal{\vec{B}}_j$ may be beneficial. In the computed examples in Section \ref{sec5}, restarting is not required since only a few $k$ steps are carried out. The {\tt nested}$\_$tGKT$_p$ method that applies the nested t-Krylov subspace method with Tikhonov regularization to solve \eqref{I6} has been described in \cite{RU1}. This method is shown in \cite{RU1} to be competitive in terms of speed and accuracy.

\vspace{.3cm}
\begin{algorithm}[H]
\SetAlgoLined
\KwIn{$\mathcal{A}, \mathcal{\vec{B}}_1, \mathcal{\vec{B}}_2, \dots, \mathcal{\vec{B}}_p$,
$\delta_1, \delta_2, \dots, \delta_p, \mathcal{L}, \eta > 1$, $k_\text{init}= 2$}
$k \leftarrow k_\text{init}$, $[\mathcal{\vec{Q}}_1, {\bf z_1}] \leftarrow
\mathtt{Normalize}(\mathcal{\vec{B}}_1)$ by Algorithm \ref{normalize}\\
Compute $ \mathcal{W}_k, \mathcal{Q}_{k+1}$ and $\mathcal{\mathcal{\bar{P}}}_k$ by
Algorithm \ref{normalize} with reorthogonalization of the tensor columns of $\mathcal{W}_k$ and $\mathcal{Q}_{k+1}$\\
Solve the minimization problem
\[
\min_{\mathcal{\vec{Z}}\in \mathbb{R}^{k \times 1 \times n}} \|\mathcal{\mathcal{\bar{P}}}_k*
\mathcal{\vec{Y}} - \mathcal{Q}_{k+1}^T* \mathcal{\vec{B}}_1 \|_F
\]
for $\mathcal{\vec{Y}}_k$ by using Algorithm \ref{backslash}\\
\While{$\|\mathcal{\mathcal{\bar{P}}}_k*\mathcal{\vec{Y}}_k -
\mathcal{Q}_{k+1}^T* \mathcal{\vec{B}}_1 \|_F \geq \eta \delta_1$}{$k \leftarrow k+1$\\
$\mathtt{Go \; to \; step \; 2}$}

Compute $\mathcal{\vec{X}}_{1,k} \leftarrow \mathcal{W}_k*\mathcal{\vec{Y}}_{1,k}$\\

\For{$j = 2,3,\dots, p$}{
$[\mathcal{\vec{Q}}_1, {\bf z_1}] \leftarrow \mathtt{Normalize}(\mathcal{\vec{B}}_j)$\\
\While{$\|\mathcal{\bar{P}}_k*\mathcal{\vec{Y}}_k - \mathcal{Q}_{k+1}^T* \mathcal{\vec{B}}_j \|_F \geq \eta \delta_j$}{
$k \leftarrow k+1$\\
\mbox{Repeat~steps~2-3 with the present tensors $\mathcal{\bar{P}}_k$,
$\mathcal{Q}_{k+1}^T$, and $\mathcal{\vec{B}}_j$}}
Compute $\mathcal{\vec{X}}_{j,k} \leftarrow \mathcal{W}_k*\mathcal{\vec{Y}}_{j,k}$\\
}
\caption{The {\tt nested}$\_$tGKB$_p$ method for the approximate solution of \eqref{I6} by solving the
problems \eqref{33no} using nested t-Krylov subspaces}\label{nestKry}
\end{algorithm}
\vspace{.3cm}

\section{Numerical examples} \label{sec5}

This section discusses the performance of the described methods. Example \ref{EX0s} illustrates the advantage of determining the solution of \eqref{I3} and \eqref{I6} by truncated iteration over Tikhonov regularization.
The purpose of the computed Example \ref{EX1s} is to discuss the performance of the T-tEVD, T-tSVD, RT-tSVD, tGKB, and t-Lanczos methods when applied to the restoration of the  gray-scale {\tt Telescope} image. Application of the T-tEVD, T-tSVD, RT-tSVD, tGKB$_p$, {\tt nested}$\_$tGKB$_p$ and t-Lanczos$_p$ methods to color image and gray-scale video restorations are considered in Examples \ref{EX2s} and \ref{EX3s}, respectively. These methods are compared to the tGKT, tGKT$_p$, and {\tt nested}$\_$tGKT$_p$ methods described in \cite{RU1} for the solution of the tensor equations \eqref{I3} and \eqref{I6}. We illustrate in Examples \ref{EX1s}-\ref{EX3s} that

\begin{enumerate}
\item accurate approximations of a few of the largest t-singular triplets of a large tensor $\mathcal{A}$ can be computed quite inexpensively by using the tGKB process. The tGKB method for \eqref{I3}, the tGKB$_p$ and {\tt nested}$\_$tGKB$_p$ methods for \eqref{I6} are seen to yield the best quality restorations and may require fewer iterations than the T-tSVD method. The latter method requiring less CPU time depends on the noise level in the data. The quality of computed solution by the tGKB, tGKB$_p$, and {\tt nested}$\_$tGKB$_p$ methods are seen to improve with Tikhonov regularization implemented by the tGKT, tGKT$_p$, and {\tt nested}$\_$tGKT$_p$ methods, respectively.

\item by using the t-Lanczos process, accurate approximations of a few of the largest t-eigenpairs of $\mathcal{A}$ can be computed fairly inexpensively. Among the methods considered, the t-Lanczos and t-Lanczos$_p$ methods are the fastest and require the least number of iterations. Both methods are implemented with reorthogonalization, and they give restorations of the worst quality. We show below that the performance of t-Lanczos-type methods can be improved tremendously by applying Tikhonov regularization. We remark that t-Lanczos-based Tikhonov regularization is yet to be considered in the literature. However, its implementation is analogous to the solution methods described in \cite{RU1, RU2}. We refer to the t-Lanczos and t-Lanczos$_p$ methods with Tikhonov regularization as t-LanczosTik and t-LanczosTik$_p$ methods, respectively.

\item randomization can be used to speed up and improve the accuracy of the tSVD method. The newly proposed RT-tSVD method is found to outperform the T-tSVD method in this regard. The implementations of the T-tSVD and RT-tSVD are similar. The T-tSVD method proceeds to compute the full tSVD of $\mathcal{A}$, then uses the discrepancy principle to determine where to truncate, whereas the RT-tSVD method computes a low tubal rank approximation of $\mathcal{A}$ by Algorithm \ref{Rtsvd} with a user-specified tolerance $\epsilon$, and determines the truncation index by the discrepancy principle.
\end{enumerate}
The factorization of $\mathcal{A}$ by the tEVD, tSVD, and R-tSVD are independent of the data tensors $\mathcal{B}$ and $\mathcal{\vec{B}}$, where
\[
\mathcal{B} = \mathcal{B}_{\rm true} + \mathcal{E},~~ {\rm and}~~ \mathcal{\vec{B}} = \mathcal{\vec{B}}_{\rm true} + \mathcal{\vec{E}}.
\]
The error tensors $\mathcal{E}$ and $\mathcal{\vec{E}}$ simulate the noise in $\mathcal{B}$ and $\mathcal{\vec{B}}$, respectively. They are defined by
\begin{equation}\label{e1}
\mathcal{E} := \widetilde{\delta} \frac{\mathcal{E}_0}{\| \mathcal{E}_0\|_F}\|\mathcal{B}_\text{true}\|_F, ~~ {\rm and}~~ \mathcal{\vec{E}} := \widetilde{\delta} \frac{\mathcal{\vec{E}}_0}{\| \mathcal{\vec{E}}_0\|_F}\|\mathcal{\vec{B}}_\text{true}\|_F,
\end{equation}
where $\widetilde{\delta}$ denotes the noise level. The entries of $\mathcal{E}_0$ and $\mathcal{\vec{E}}_0$ are drawn from a standard normal distribution. The effectiveness of the solution methods is determined for two noise levels: $\widetilde{\delta} = 10^{-2}$ and $\widetilde{\delta} = 10^{-3}$. 

\begin{figure}[!htb]
\hspace{-.8cm}
\minipage{0.37\textwidth}
\includegraphics[width=\linewidth]{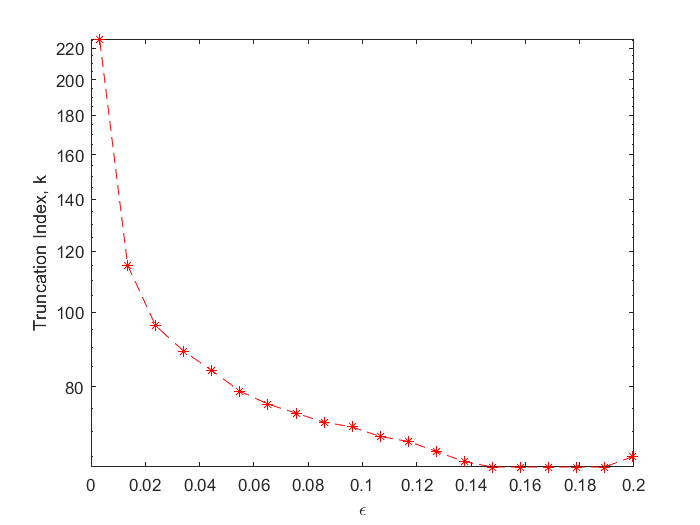} 
\endminipage\hfill \hspace{-1.3cm}
\minipage{0.37\textwidth}
\includegraphics[width=\linewidth]{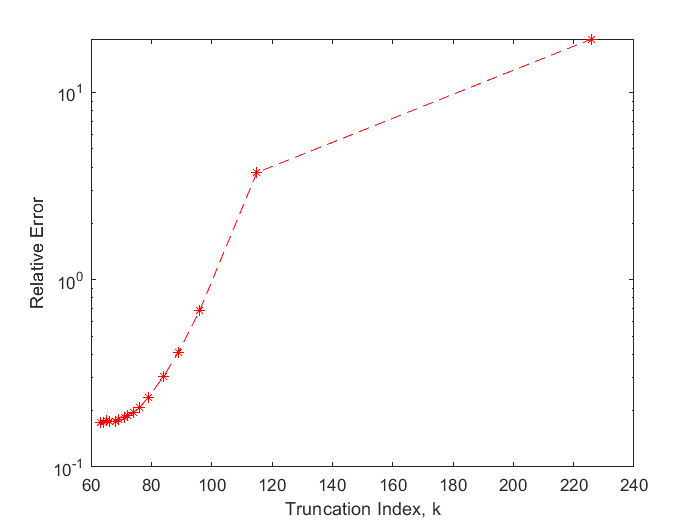}
\endminipage\hfill \hspace{-1.3cm}
\minipage{0.37\textwidth}
\includegraphics[width=\linewidth]{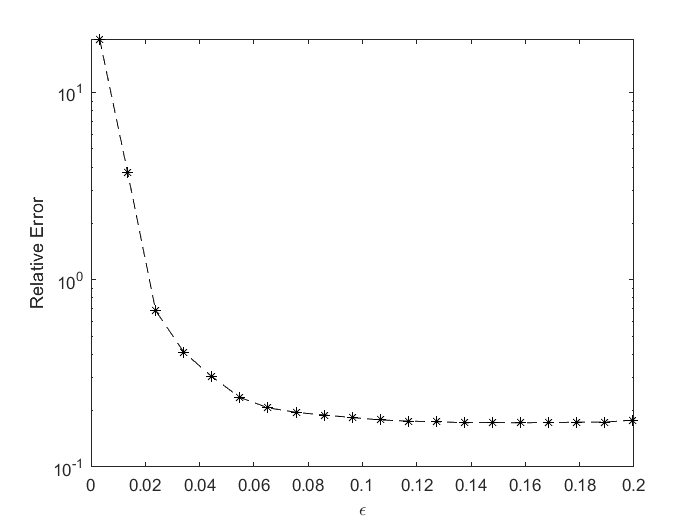}
\endminipage \hspace{-.8cm} \vspace{-.3cm}
\caption{\small Sensistivity of the RT-tSVD method to $\epsilon$ and $k$ when applied to the restoration of {\tt Telescope} image for $\widetilde{\delta}=10^{-2}$.}
\label{Fig: 6sT1}
\end{figure}

\begin{figure}[!htb]
\hspace{-.8cm}
\minipage{0.37\textwidth}
\includegraphics[width=\linewidth]{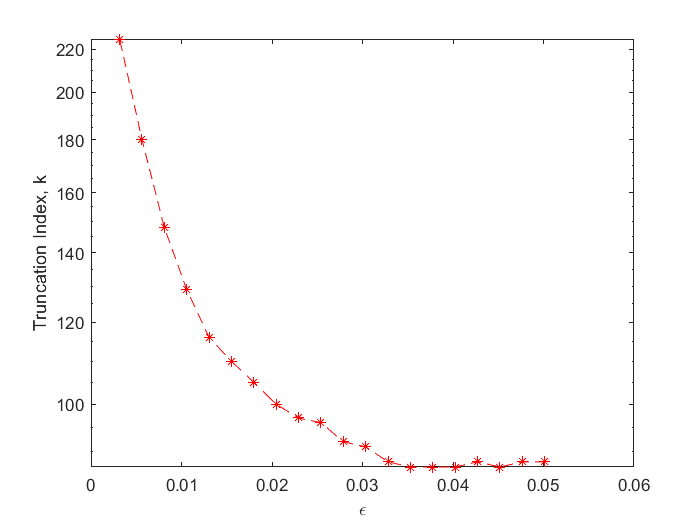} 
\endminipage\hfill \hspace{-1.3cm}
\minipage{0.37\textwidth}
\includegraphics[width=\linewidth]{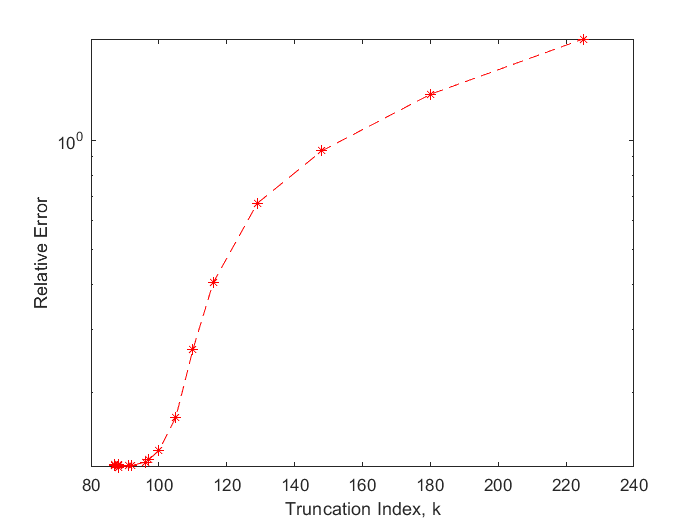}
\endminipage\hfill \hspace{-1.3cm}
\minipage{0.37\textwidth}
\includegraphics[width=\linewidth]{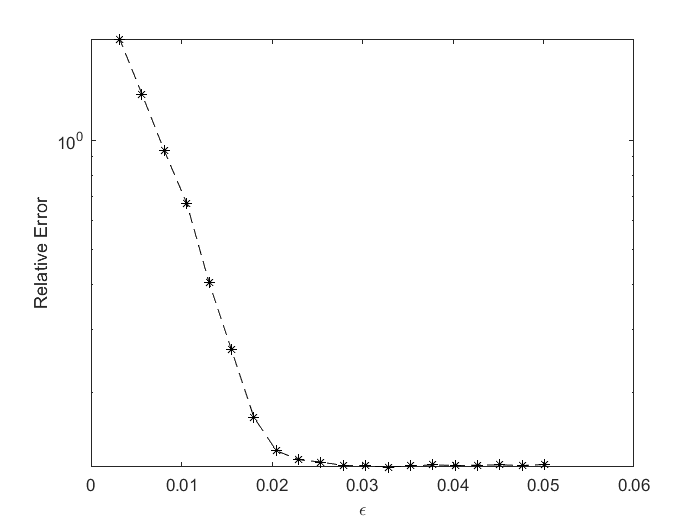}
\endminipage \hspace{-.8cm} \vspace{-.3cm}
\caption{\small Sensistivity of the RT-tSVD method to $\epsilon$ and $k$ when applied to the restoration of {\tt Telescope} image for $\widetilde{\delta}=10^{-3}$.}
\label{Fig: 6sT01}
\end{figure}

\begin{figure}[!htb]
\hspace{-.8cm}
\minipage{0.37\textwidth}
\includegraphics[width=\linewidth]{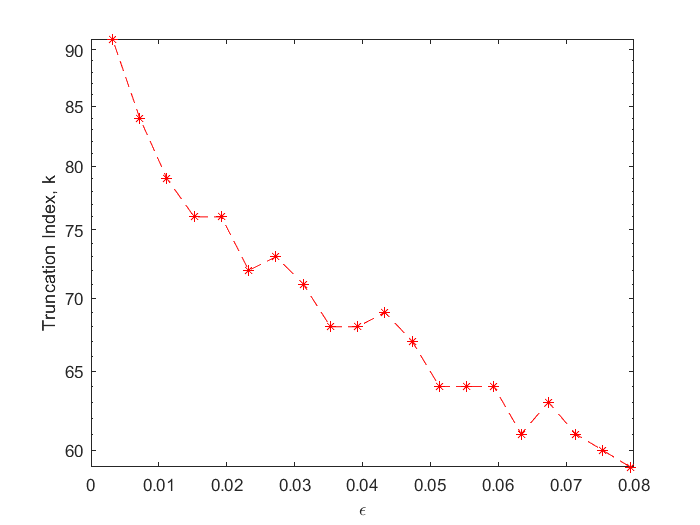} 
\endminipage\hfill \hspace{-1.3cm}
\minipage{0.37\textwidth}
\includegraphics[width=\linewidth]{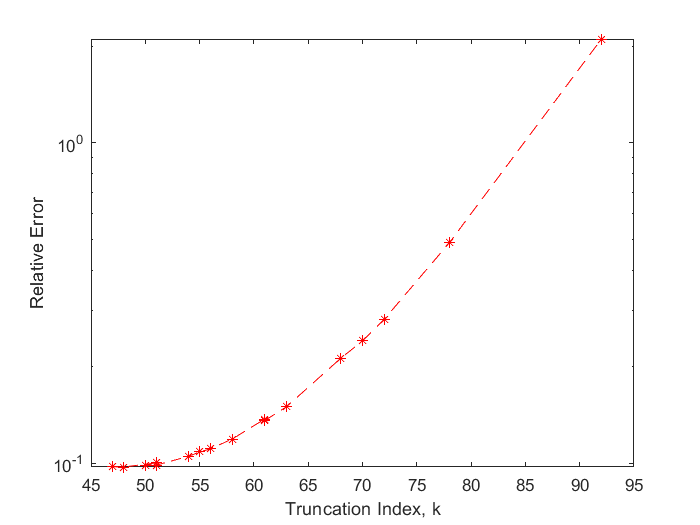}
\endminipage\hfill \hspace{-1.3cm}
\minipage{0.37\textwidth}
\includegraphics[width=\linewidth]{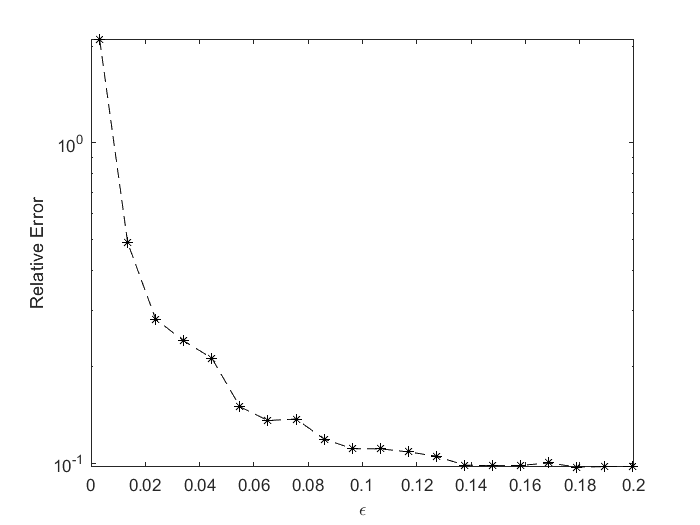}
\endminipage \hspace{-.8cm} \vspace{-.3cm}
\caption{\small Sensistivity of the RT-tSVD method to $\epsilon$ and $k$ when applied to the restoration of {\tt papav256} image for $\widetilde{\delta}=10^{-2}$.}
\label{Fig: 6sp1}
\end{figure}

\begin{figure}[!htb]
\hspace{-.8cm}
\minipage{0.37\textwidth}
\includegraphics[width=\linewidth]{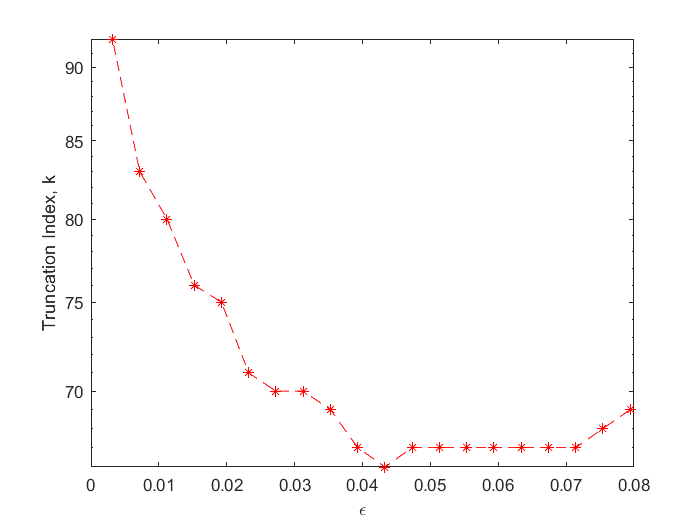} 
\endminipage\hfill \hspace{-1.3cm}
\minipage{0.37\textwidth}
\includegraphics[width=\linewidth]{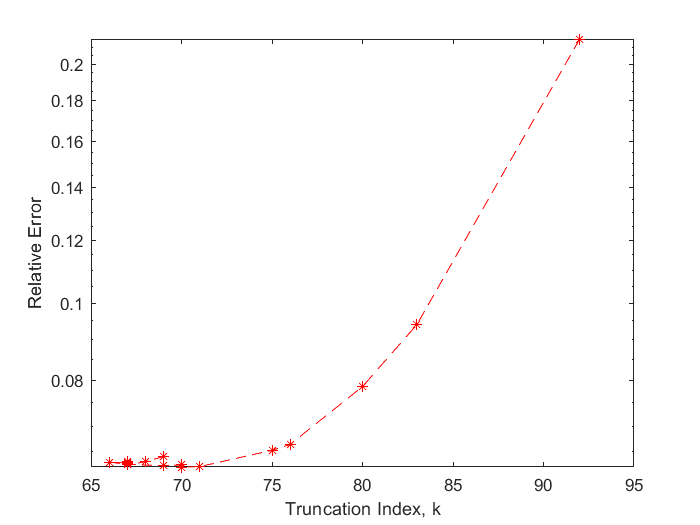}
\endminipage\hfill \hspace{-1.3cm}
\minipage{0.37\textwidth}
\includegraphics[width=\linewidth]{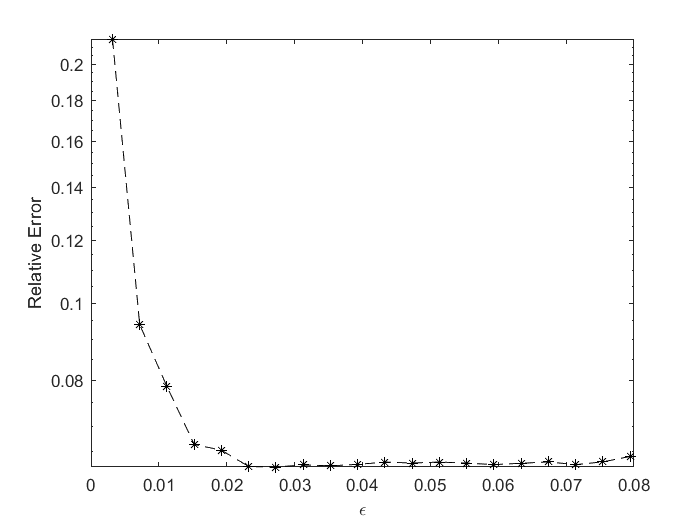}
\endminipage \hspace{-.8cm} \vspace{-.3cm}
\caption{\small Sensistivity of the RT-tSVD method to $\epsilon$ and $k$ when applied to the restoration of {\tt papav256} image for $\widetilde{\delta}=10^{-3}$.}
\label{Fig: 6sp2}
\end{figure}

\begin{figure}[!htb]
\hspace{-.8cm}
\minipage{0.37\textwidth}
\includegraphics[width=\linewidth]{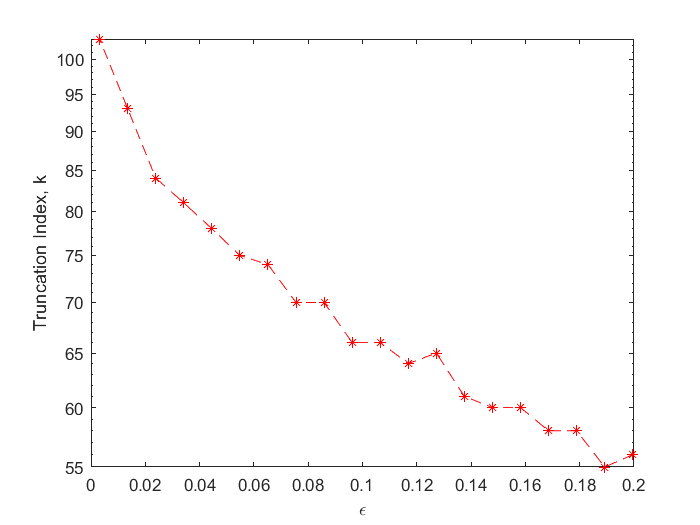} 
\endminipage\hfill \hspace{-1.3cm}
\minipage{0.37\textwidth}
\includegraphics[width=\linewidth]{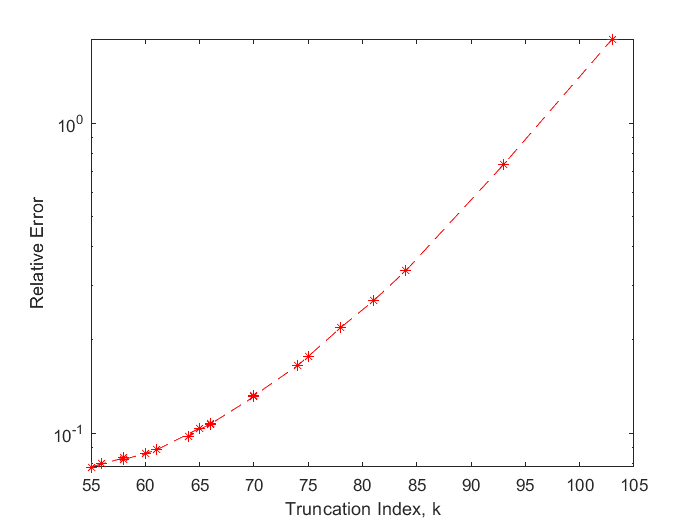}
\endminipage\hfill \hspace{-1.3cm}
\minipage{0.37\textwidth}
\includegraphics[width=\linewidth]{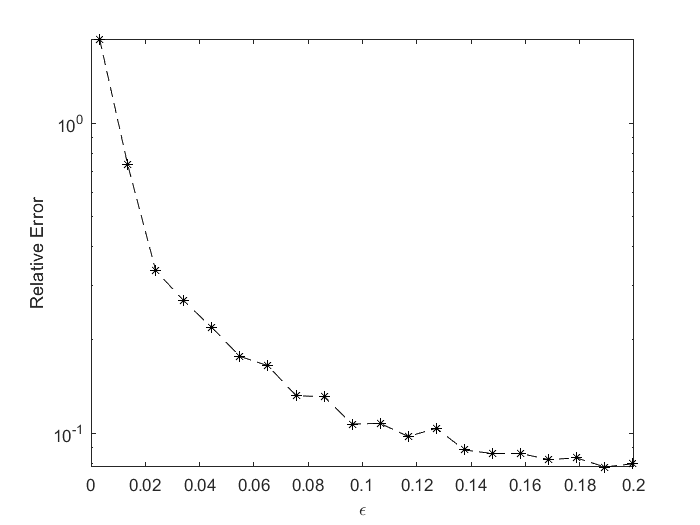}
\endminipage \hspace{-.8cm} \vspace{-.3cm}
\caption{\small Sensistivity of the RT-tSVD method to $\epsilon$ and $k$ when applied to the restoration of gray-scale {\tt video} for $\widetilde{\delta}=10^{-2}$.}
\label{Fig: 6sv1}
\end{figure}

\begin{figure}[!htb]
\hspace{-.8cm}
\minipage{0.37\textwidth}
\includegraphics[width=\linewidth]{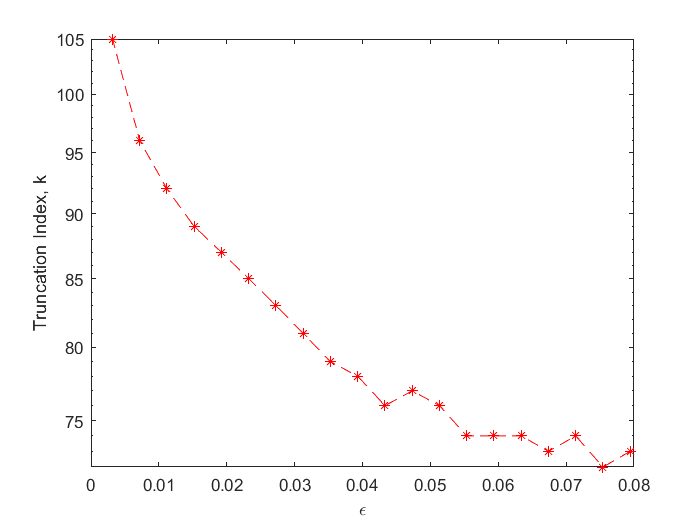} 
\endminipage\hfill \hspace{-1.3cm}
\minipage{0.37\textwidth}
\includegraphics[width=\linewidth]{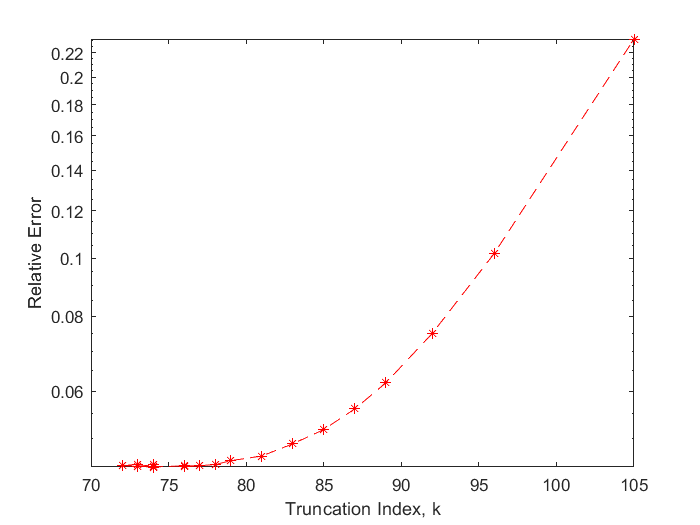}
\endminipage\hfill \hspace{-1.3cm}
\minipage{0.37\textwidth}
\includegraphics[width=\linewidth]{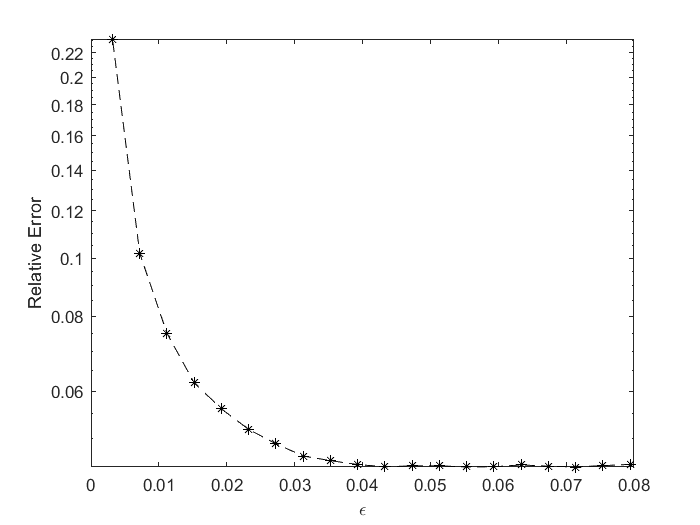}
\endminipage \hspace{-.8cm} \vspace{-.3cm}
\caption{\small Sensistivity of the RT-tSVD method to $\epsilon$ and $k$ when applied to the restoration of gray-scale {\tt video} for $\widetilde{\delta}=10^{-3}$.}
\label{Fig: 6sv2}
\end{figure}

The choice of $\epsilon$ in Algorithm \ref{Rtsvd} for the RT-tSVD method depends on the level of noise in the data. A smaller noise level often correspond to smaller value of $\epsilon$. This ensures that the low tubal rank approximation of $\mathcal{A}$ determined by Algorithm \ref{Rtsvd} does not overestimate the sizes of the factors $\mathcal{Q}$ and $\mathcal{\tilde{B}}$. We see from Figures \ref{Fig: 6sT1} - \ref{Fig: 6sv2} that the smaller the $\epsilon$, the larger the truncation index, $k$, determined by the discrepancy principle. Moreover, very large $k$ values result in low quality reconstructions determined by the relative error. Generally, the relative error decays polynomially with increasing $\epsilon$  values. Figures \ref{Fig: 6sT1}-\ref{Fig: 6sv2} inform our choices of $\epsilon$ below. 

Let $r$ be the number of iterations required by Algorithm \ref{Rtsvd} to satisfy the inequality \eqref{apperr}, and let $\rho$ be the oversampling parameter. The initial truncation index $k = k_{\tt init}$ for the RT-tSVD method is determined by
\be \label{kinit}
k_{\tt init} = r - \rho,
\ee
where $k_{\tt init}$ is increased until the discrepancy principle \eqref{disp} is satisfied. We take $\rho = 3$ in Example \ref{EX0s}, and $\rho = 10$ in Examples \ref{EX1s}-\ref{EX3s}. The choices of $\epsilon$ determined after repeated numerical experiments, and the computed values of $r$ by Algorithm \ref{Rtsvd} for Examples \ref{EX0s}-\ref{EX3s} are summarized in Table \ref{Tab: 0} below. 

\begin{table}[h!]
\begin{center}
\begin{tabular}{cccccccc}
\cmidrule(lr){1-4}
\multicolumn{1}{c}{Noise level} & \multicolumn{1}{c}{Example}&\multicolumn{1}{c}{$\epsilon$} & \multicolumn{1}{c}{$r$}
\\ \cmidrule(lr){1-4} 
\multirow{4}{2em}{$10^{-3}$} &\ref{EX0s} &$10^{-1.5}$ &3\\
&\ref{EX1s} &$10^{-1.5}$ &99\\
&\ref{EX2s} & $10^{-1.5}$ &79\\
&\ref{EX3s}& $10^{-1.2}$&81 \\ \cmidrule(lr){1-4}
\multirow{4}{2em}{$10^{-2}$} &\ref{EX0s}&$10^{-1.5}$&3 \\
&\ref{EX1s}&$10^{-0.7}$&66 \\
&\ref{EX2s}& $10^{-0.7}$&57 \\
&\ref{EX3s}& $10^{-0.7}$&66 \\ \cmidrule(lr){1-4}
\end{tabular}
\end{center} \vspace{-.5cm}
\caption{\small Parameters used by the RT-tSVD method. The number of iteration $r$ is determined by Algorithm \ref{Rtsvd} for a given tolerance $\epsilon$.}
\label{Tab: 0}
\end{table}

The relative errors for \eqref{I3} and \eqref{I6} are computed as
\begin{equation*}
E_{\rm method} = \frac{\|\mathcal{\vec{X}}_{\rm method} - \mathcal{\vec{X}}_\text{true}\|_F}{\|\mathcal{\vec{X}}_\text{true}\|_F},~~{\rm and}~~E_{\rm method} = \frac{\|\mathcal{X}_{\rm method} - \mathcal{X}_\text{true}\|_F}{\|\mathcal{X}_\text{true}\|_F},
\end{equation*}
where $\mathcal{\vec{X}}_\text{true}$ and $\mathcal{X}_\text{true}$ are the desired true solution of the error-free problems corresponding to \eqref{I3} and \eqref{I6}, respectively.

The blurring operator $\mathcal{A} \in \mathbb{R}^{\ell \times m\times n}$ with $\ell = m$ has frontal slices $\mathcal{A}^{(i)} \in \mathbb{R}^{\ell \times m}$, $i = 1,2,\dots,n$, that are generated by using the function $\mathtt{blur}$ from Hansen's Regularization Tools \cite{Haa}. Specifically,
\begin{equation*}\label{z}
z = \mathtt{[exp(-([0:band-1].^2)/(2\sigma^2)),zeros(1,N-band)]},
\end{equation*}
\begin{equation}\label{za}
A = \frac{1}{\sigma\sqrt{2\pi}} \mathtt{toeplitz}(z), \;\; \mathcal{A}^{(i)} = A(i,1)A,
~~ i = 1,2, \dots, n.
\end{equation}

The discrepancy principle \eqref{disp} is used to determine the regularization parameter (truncation index) for all methods with $\eta = 1.1$ in Examples \ref{EX0s} and \ref{EX1s}, and $\eta = 1.2$ in Examples \ref{EX2s} and \ref{EX3s}. We remark that the theory behind the discrepancy principle does not prescribe how $\eta$ should be chosen.	The Tikhonov regularization parameter $\mu$ is determined by using the bisection method over the interval $[10^{-3},10^{5}]$ in Example \ref{EX0s}, and $[10^1, 10^7]$ in Examples \ref{EX1s}-\ref{EX3s}.

\begin{Ex}\label{EX0s}
This example considers the solution of the minimization problems \eqref{I3} and \eqref{I6} with nonsymmetric tensors $\mathcal{A} \in \mathbb{R}^{500 \times 500 \times 500}$ and $\mathcal{A} \in \mathbb{R}^{300 \times 300 \times 300}$, respectively. The tensors for both problems are generated by their frontal slices as described below. We illustrate that the tGKB method for \eqref{I3}, the tGKB$_p$, and {\tt nested}$\_$tGKB$_p$ methods for \eqref{I6} by truncated iterations can yield more accurate approximate solutions than their counterpart based on Tikhonov regularization. Let $A_1 = {\tt gallery('prolate',n,0.46)}$, and $A_2 = {\tt baart}(n)$. Then
\[
\mathcal{A}^{(i)} = A_1(i,1)A_2, ~~i=1,2,\dots,n, ~~n \in \{300,500\},
\]
where $A_1$ is an ill-conditioned symmetric positive definite Toeplitz matrix from MATLAB, $A_2$ is a nonsymmteric matrix from \cite{Haa}, and $\mathcal{A}$ are generated by folding the first block column of $A_1\otimes A_2$, where $\otimes$ denotes the Kronecker product. The condition number of each nonsymmteric matrix $\mathcal{A}^{(i)}$, computed using MATLAB's {\tt cond} operator, is larger than $10^{17}$. Hence each $\mathcal{A}^{(i)}$ is numerically singular.

\begin{table}[h!]
\begin{center}
\begin{tabular}{cccccccc}
\cmidrule(lr){1-5}
\multicolumn{1}{c}{Noise level} & \multicolumn{1}{c}{Method}&\multicolumn{1}{c}{$k$} & \multicolumn{1}{c}{Relative error}& \multicolumn{1}{c}{CPU time (secs)}
\\ \cmidrule(lr){1-5}
\multirow{4}{4em}{$10^{-3}$} &tGKB &3 &$5.9830\cdot 10^{-3}$ &12.04 \\
&T-tSVD &3 & $6.0031\cdot 10^{-3}$ & 34.14\\
&RT-tSVD&3 & $5.5868\cdot 10^{-3}$ & 24.16 \\
&tGKT &3 &$1.3826\cdot 10^{-2}$ &12.51\\ \cmidrule(lr){1-5}
\multirow{4}{4em}{$10^{-2}$} &tGKB &2 &$7.1518\cdot 10^{-2}$ &4.81 \\
&T-tSVD &2 & $7.2474\cdot 10^{-2}$ & 30.36 \\
&RT-tSVD&2 & $7.2472\cdot 10^{-2}$ & 18.65 \\
&tGKT &2 &$7.3710\cdot 10^{-2}$ &5.34 \\ \cmidrule(lr){1-5}
\end{tabular}
\end{center} \vspace{-.5cm}
\caption{\small Results for tGKB, T-tSVD, RT-tSVD and tGKT methods when applied to a problem defined by {\tt baart} \cite{Haa} and prolate matrix from MATLAB.}
\label{Tab: 1s0}
\end{table}

Let $\mathcal{\vec{X}}_{\rm true} \in \mathbb{R}^{500 \times 1 \times 500}$ and $\mathcal{X}_{\rm true} \in \mathbb{R}^{300 \times 3 \times 300}$ with unit entries be the true solutions of \eqref{I3} and \eqref{I6}, respectively, and determine the associated noise-contaminated right-hand side by $\mathcal{\vec{B}} = \mathcal{A}*\mathcal{\vec{X}}_\text{true} + \mathcal{\vec{E}}$, and $\mathcal{B} = \mathcal{A}*\mathcal{X}_\text{true} + \mathcal{E}$, respectively, where the noise tensors $\mathcal{\vec{E}}$ and $\mathcal{E}$ are defined above.

Table \ref{Tab: 1s0} displays the number of iterations, relative error, and CPU time determined by the tGKB, T-tSVD, RT-tSVD, and tGKT methods for the approximate solution of \eqref{I3}. The tGKT method that is based on Tikhonov regularization determines approximate solutions of the worst quality for both noise levels. This method computes the regularization parameters $\mu = 1.02 \cdot 10^3$ and $\mu = 78$ for $\widetilde{\delta} = 10^{-3}$ and $\widetilde{\delta} = 10^{-2}$, respectively. For both noise levels, the tGKT and tGKB methods that use the bidiagonalization process are the fastest. While the tGKB method determines an approximate solution of the best quality for $\widetilde{\delta} = 10^{-2}$, the RT-tSVD gives the approximate solution of higher quality for $\widetilde{\delta} = 10^{-3}$. The RT-tSVD method determines solution of highest quality than the T-tSVD method for both noise levels.

Similarly, Table \ref{Tab: 3s0} compares the tGKB$_p$, {\tt nested}$\_$tGKB$_p$, T-tSVD, RT-tSVD, {\tt nested}$\_$tGKT$_p$ and tGKT$_p$ methods for the approximate solution of \eqref{I6}. Here and below, the table entries $k = (k_1, k_2, \dots, k_j)$, $j=1,2,\dots,p$, indicates that the chosen $p$-method computes $k$ different or the same number of iterations corresponding to each lateral slice $\mathcal{\vec{B}}_j$ of the data $\mathcal{B}$. Additionally, the computed Tikhonov regularization parameters determined by the $p$-methods, each corresponding to $\mathcal{\vec{B}}_j$, are omitted for convenience since each of the $p$-methods computes $p$ different regularization parameters.

Table \ref{Tab: 3s0} shows that Tikhonov regularization-based tGKT$_p$ and {\tt nested}$\_$tGKT$_p$ methods determine solution of the worst quality for both noise levels. The {\tt nested}$\_$tGKB$_p$ method that uses nested t-Krylov subspaces is the fastest for both noise levels. The RT-tSVD and tGKB$_p$ methods give approximate solutions of highest quality for $\widetilde{\delta} = 10^{-3}$ and $\widetilde{\delta} = 10^{-2}$, respectively. The tGKT$_p$ and T-tSVD methods are the slowest for $\widetilde{\delta} = 10^{-3}$ and $\widetilde{\delta} = 10^{-2}$, respectively.

\begin{table}[h!]
\begin{center}
\begin{tabular}{cccccccc}
\cmidrule(lr){1-5}
\multicolumn{1}{c}{Noise level} & \multicolumn{1}{c}{Method}&\multicolumn{1}{c}{$k$} & \multicolumn{1}{c}{Relative error}& \multicolumn{1}{c}{CPU time (secs)}
\\ \cmidrule(lr){1-5}
\multirow{6}{4em}{$10^{-3}$} &tGKB$_p$ &(3,3,3)&$6.1528\cdot 10^{-3}$ &8.01 \\
&{\tt nested}$\_$tGKB$_p$ &3 & $6.1544\cdot 10^{-3}$ & 2.79\\
&T-tSVD &3 & $6.1617\cdot 10^{-3}$ & 6.53 \\
&RT-tSVD&3 & $5.9258\cdot 10^{-3}$ & 5.46\\
&tGKT$_p$&(3,3,3) &$1.3989\cdot 10^{-2}$ &9.12 \\
&{\tt nested}$\_$tGKT$_p$ &3& $2.3087\cdot 10^{-2}$ & 5.66\\ \cmidrule(lr){1-5}
\multirow{6}{4em}{$10^{-2}$} &tGKB$_p$ &(2,2,2) &$7.1541\cdot 10^{-2}$ &3.18 \\
&{\tt nested}$\_$tGKB$_p$ &2 & $7.1547\cdot 10^{-2}$ & 1.14\\
&T-tSVD &2 & $7.2494\cdot 10^{-2}$ & 5.31\\
&RT-tSVD &2 & $7.2481\cdot 10^{-2}$ & 4.56 \\
&tGKT$_p$&(2,2,2)& $7.3753\cdot 10^{-2}$ & 4.19 \\
&{\tt nested}$\_$tGKT$_p$ &2& $8.1115\cdot 10^{-2}$ & 3.85\\ \cmidrule(lr){1-5}
\end{tabular}
\end{center} \vspace{-.5cm}
\caption{\small Results for tGKB$_p$, {\tt nested}$\_$tGKB$_p$, T-tSVD, RT-tSVD, tGKT$_p$ and {\tt nested}$\_$tGKT$_p$ methods when applied to a problem defined by {\tt baart} \cite{Haa} and prolate matrix from MATLAB.}
\label{Tab: 3s0}
\end{table}
\end{Ex}

\begin{Ex} (2-D image restoration)\label{EX1s} This example considers a two-dimensional image restoration problem. Here, the true {\tt Telescope} image is shown on the left-hand side of Figure \ref{Fig: 1s}. This image is stored as a tensor slice $\mathcal{\vec{X}}_{\rm true}\in \mathbb{R}^{256 \times 1 \times 256}$ and blurred by $\mathcal{A}$ in \eqref{za} with $\sigma=3$ and ${\tt band} = 9$. The blurred and noisy image, shown in Figure \ref{Fig: 1s} (right) for a noise level $\widetilde{\delta} = 10^{-3}$, is generated by $\mathcal{\vec{B}} = \mathcal{A}*\mathcal{\vec{X}}_\text{true} + \mathcal{\vec{E}}$ where $\mathcal{\vec{E}}$ is defined by the right-hand side of \eqref{e1}.

The relative error, CPU time and number of iterations required by each method are presented in Table \ref{Tab: 1s} for both noise levels. The restored images by the tGKB, t-Lanczos and RT-tSVD methods are shown in Figure \ref{Fig: 2s} (left-right) for $\widetilde{\delta} = 10^{-3}$. Table \ref{Tab: 1s} shows that the tGKB and tGKT methods are the slowest for $\widetilde{\delta} = 10^{-3}$, but they are faster than the T-tEVD, T-tSVD and RT-tSVD methods for $\widetilde{\delta} = 10^{-2}$. The t-Lanczos and t-LanczosTik methods are the fastest, but the t-Lanczos method yields restoration of the worst quality among all methods for both noise levels. The RT-tSVD and t-Lanczos-type methods require the most and the least number of iterations for both noise levels, respectively. The RT-tSVD method is faster than the T-tSVD method for both noise levels. While the RT-tSVD method yields restoration of higher accuracy than the T-tSVD method for $\widetilde{\delta} = 10^{-3}$, the latter method gives more quality reconstruction for $\widetilde{\delta} = 10^{-2}$.

\begin{figure}[htb]
\hspace{1cm}
\minipage{0.42\textwidth}
\includegraphics[width=\linewidth]{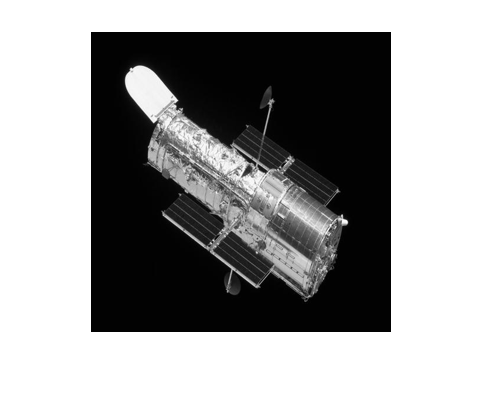} 
\endminipage\hfill \hspace{-2cm}
\minipage{0.42\textwidth}
\includegraphics[width=\linewidth]{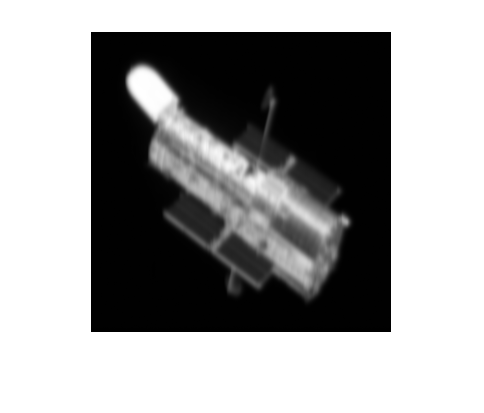}
\endminipage\hfill \hspace{-1cm} \vspace{-.9cm}
\caption{\small True {\tt Telescope} image (left), and blurred and noisy {\tt Telescope} image (right) for $\widetilde{\delta}=10^{-3}$.}
\label{Fig: 1s}
\end{figure}

\begin{figure}[htb]
\hspace{-1cm}
\minipage{0.42\textwidth}
\includegraphics[width=\linewidth]{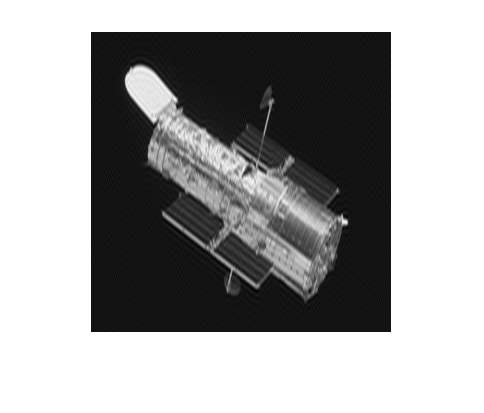} 
\endminipage\hfill \hspace{-2cm}
\minipage{0.42\textwidth}
\includegraphics[width=\linewidth]{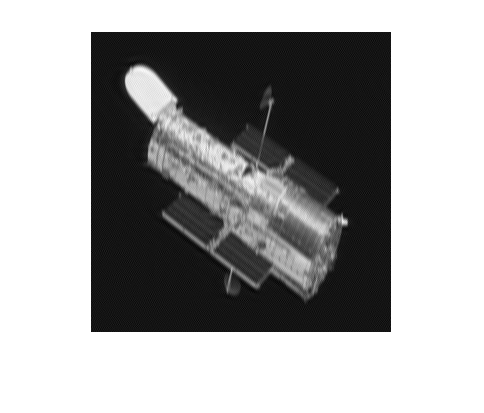}
\endminipage\hfill \hspace{-2cm}
\minipage{0.42\textwidth}
\includegraphics[width=\linewidth]{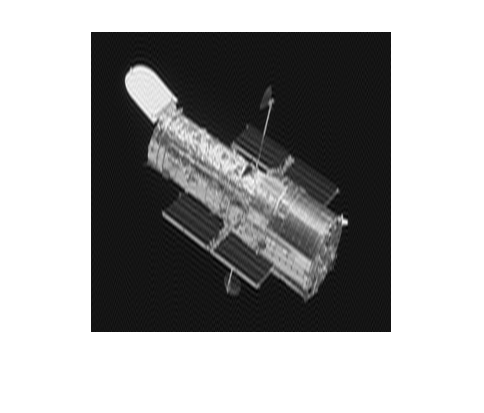}
\endminipage \hspace{-1cm} \vspace{-.9cm}
\caption{\small Reconstructed {\tt Telescope} images by tGKB method after $35$ iterations
(left), t-Lanczos method after $6$ iterations (middle), and RT-tSVD method after $89$ iterations
(right) for $\widetilde{\delta}=10^{-3}$.}
\label{Fig: 2s}
\end{figure}

\begin{table}[h!]
\begin{center}
\begin{tabular}{cccccccc}
\cmidrule(lr){1-5}
\multicolumn{1}{c}{Noise level} & \multicolumn{1}{c}{Method}&\multicolumn{1}{c}{$k$} & \multicolumn{1}{c}{Relative error}& \multicolumn{1}{c}{CPU time (secs)}
\\ \cmidrule(lr){1-5}
\multirow{7}{4em}{$10^{-3}$} &tGKB &35 &$1.2235\cdot 10^{-1}$ &377.33 \\
&t-Lanczos &6 & $1.5933\cdot 10^{-1}$ &5.77 \\
&T-tEVD&86& $1.2575\cdot 10^{-1}$ & 117.68 \\
&T-tSVD &86 & $1.2575\cdot 10^{-1}$ & 103.48\\
&RT-tSVD&89 & $1.2475\cdot 10^{-1}$ & 67.92 \\
&tGKT &35 &$1.1975\cdot 10^{-1}$ &434.55 \\
&t-LanczosTik &6 & $1.4555\cdot 10^{-1}$ &6.94 \\ \cmidrule(lr){1-5}
\multirow{7}{4em}{$10^{-2}$} &tGKB &7 &$1.6953\cdot 10^{-1}$ &16.89 \\
&t-Lanczos &3 & $4.5149\cdot 10^{-1}$ &1.50 \\
&T-tEVD&62& $1.7150\cdot 10^{-1}$ & 87.66 \\
&T-tSVD &62 & $1.7150\cdot 10^{-1}$ & 73.33 \\
&RT-tSVD&64 & $1.7622\cdot 10^{-1}$ & 53.14 \\
&tGKT &7 &$1.5394\cdot 10^{-1}$ &18.16 \\
&t-LanczosTik &3& $1.6732\cdot 10^{-1}$ &1.92\\ \cmidrule(lr){1-5}
\end{tabular}
\end{center} \vspace{-.5cm}
\caption{\small Results for tGKB, t-Lanczos, T-tEVD, T-tSVD, RT-tSVD, tGKT, and t-LanczosTik methods when applied to the restoration of {\tt Telescope} image.}
\label{Tab: 1s}
\end{table}
\end{Ex}

\begin{Ex}{(Color image deblurring)}\label{EX2s}
This example discusses color image restoration. The frontal slices $\mathcal{A}^{(i)}$, $i = 1,2,\dots,n$, of $\mathcal{A} \in \mathbb{R}^{256 \times 256 \times 256}$ are generated by \eqref{za} with {\tt band} = 12 and $\sigma=3$. The true unknown image, shown on the left-hand side of Figure \ref{Fig: 3s}, is the three-channel $\mathtt{papav256}$ RGB image. This image is stored as a tensor $\mathcal{X}_{\text{true}} \in \mathbb{R}^{256 \times 3 \times 256}$, and blurred by $\mathcal{A}$ above. The blurred and noisy image represented by $\mathcal{B} \in \mathbb{R}^{256 \times 3 \times 256}$ is generated by $\mathcal{B} = \mathcal{A}*\mathcal{X}_\text{true} + \mathcal{E}$, and displayed on the right-hand side of Figure \ref{Fig: 3s} for $\widetilde{\delta} = 10^{-3}$, where $\mathcal{E}$ is generated as described by the left-hand side of \eqref{e1}. The restored images by the tGKB$_p$, {\tt nested}$\_$tGKB$_p$, and t-Lanczos$_p$ methods for noise level of $\widetilde{\delta} = 10^{-3}$ are shown in Figure \ref{Fig: 4s}.

Table \ref{Tab: 2s} shows the number of iterations, relative error, and the CPU time required for each method. For the truncated iteration methods, the tGKB$_p$ method yields restorations of the highest quality for $\widetilde{\delta} = 10^{-3}$, while the T-tSVD and T-tEVD methods give the most accurate reconstructions for $\widetilde{\delta} = 10^{-2}$. The t-Lanczos$_p$ method is the fastest and it yields restoration of the worst quality for both noise levels. The quality of the reconstruction by the RT-tSVD method is more accurate than those of the T-tSVD method for $\widetilde{\delta} = 10^{-3}$. The RT-tSVD method is faster than the T-tSVD method for both noise levels. The t-Krylov recycling approaches implemented by the {\tt nested}$\_$tGKB$_p$ and {\tt nested}$\_$tGKT$_p$ methods are quite competitive in terms of speed and quality of reconstructions. Both methods reduce the CPU time required by the tGKB$_p$ and tGKT$_p$ methods, respectively, by more than half for $\widetilde{\delta} = 10^{-3}$.

\begin{figure}[!htb]
\hspace{1cm}
\minipage{0.42\textwidth}
\includegraphics[width=\linewidth]{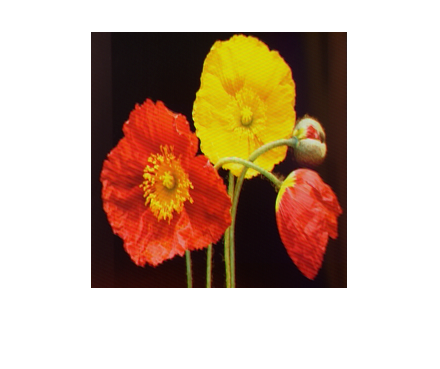} 
\endminipage\hfill \hspace{-2cm}
\minipage{0.42\textwidth}
\includegraphics[width=\linewidth]{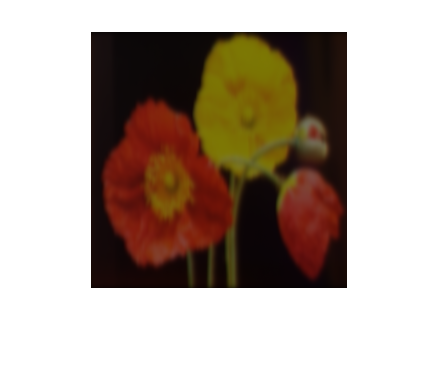}
\endminipage\hfill \hspace{-1cm} \vspace{-.9cm}
\caption{\small True {\tt papav256} image (left), and blurred and noisy {\tt papav256} image (right) for $\widetilde{\delta}=10^{-3}$.}
\label{Fig: 3s}
\end{figure}

\begin{figure}[!htb]
\hspace{-1cm}
\minipage{0.42\textwidth}
\includegraphics[width=\linewidth]{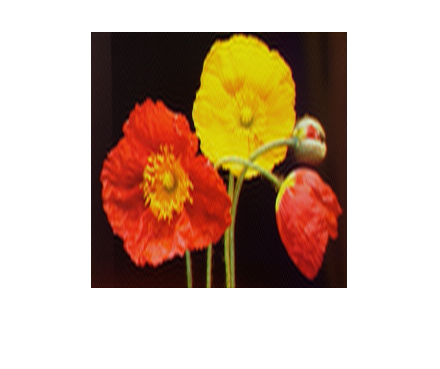} 
\endminipage\hfill \hspace{-2cm}
\minipage{0.42\textwidth}
\includegraphics[width=\linewidth]{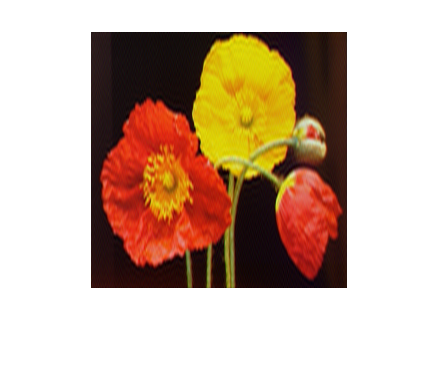}
\endminipage\hfill \hspace{-2cm}
\minipage{0.42\textwidth}
\includegraphics[width=\linewidth]{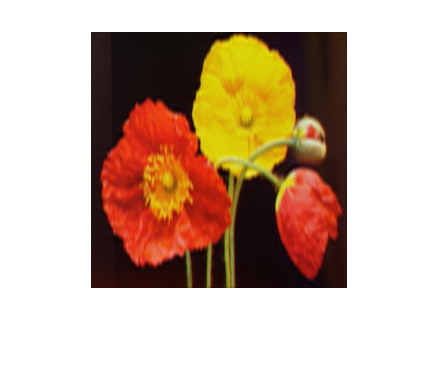}
\endminipage \hspace{-1cm} \vspace{-.9cm}
\caption{\small Reconstructed {\tt papav256} images by tGKB$_p$ method (left), {\tt nested}$\_$tGKB$_p$ method after $23$ iterations (middle), and t-Lanczos$_p$ method (right) for $\widetilde{\delta}=10^{-3}$.}
\label{Fig: 4s}
\end{figure}

\begin{table}[h!]
\begin{center}
\begin{tabular}{cccccccc}
\cmidrule(lr){1-5}
\multicolumn{1}{c}{Noise level} & \multicolumn{1}{c}{Method}&\multicolumn{1}{c}{$k$} & \multicolumn{1}{c}{Relative error}& \multicolumn{1}{c}{CPU time (secs)}
\\ \cmidrule(lr){1-5}
\multirow{9}{4em}{$10^{-3}$} &tGKB$_p$ &(20,22,19)&$6.1484\cdot 10^{-2}$ &243.21 \\
&{\tt nested}$\_$tGKB$_p$ &23 & $6.1618\cdot 10^{-2}$ &110.88 \\
&t-Lanczos$_p$ &(5,5,5) & $8.9132\cdot 10^{-2}$ &7.77\\
&T-tEVD&66& $6.2908\cdot 10^{-2}$ & 60.89 \\
&T-tSVD &66 & $6.2908\cdot 10^{-2}$ & 52.50 \\
&RT-tSVD&69 & $6.2131\cdot 10^{-2}$ & 34.14 \\
&tGKT$_p$&(20,22,19)& $5.9602\cdot 10^{-2}$ & 263.75 \\
&{\tt nested}$\_$tGKT$_p$ &23& $5.9557\cdot 10^{-2}$ & 125.49 \\
&t-LanczosTik$_p$&(5,5,5) & $6.7251\cdot 10^{-2}$ & 9.36 \\ \cmidrule(lr){1-5}
\multirow{9}{4em}{$10^{-2}$} &tGKB$_p$ &(5,4,4) &$9.8184\cdot 10^{-2}$ &12.57 \\
&{\tt nested}$\_$tGKB$_p$ &6 & $1.3160\cdot 10^{-1}$ &8.32 \\
&t-Lanczos$_p$ &(2,2,2) & $2.7783\cdot 10^{-1}$ &1.18 \\
&T-tEVD&45& $9.4964\cdot 10^{-2}$ & 41.58 \\
&T-tSVD &45 & $9.4964\cdot 10^{-2}$ & 34.39 \\
&RT-tSVD&47 & $9.7710\cdot 10^{-2}$ & 24.96 \\
&tGKT$_p$&(5,4,4)& $8.3342\cdot 10^{-2}$ & 15.22 \\
&{\tt nested}$\_$tGKT$_p$ &6& $1.1984\cdot 10^{-1}$ & 12.88 \\
&t-LanczosTik$_p$&(2,2,2) & $1.0019\cdot 10^{-1}$ & 2.30 \\ \cmidrule(lr){1-5}
\end{tabular}
\end{center} \vspace{-.5cm}
\caption{\small Results for tGKB$_p$, {\tt nested}$\_$tGKB$_p$, t-Lanczos$_p$, T-tEVD, T-tSVD, RT-tSVD, tGKT$_p$, {\tt nested}$\_$tGKT$_p$ and t-LanczosTik$_p$ methods when applied to the restoration of {\tt papav256} image.}
\label{Tab: 2s}
\end{table}

\end{Ex}

\begin{Ex}{(Gray-scale video restoration)}\label{EX3s} This example illustrates the performance of the described methods when applied to the restoration of six consecutive gray-scale {\tt video} frames from MATLAB. Each frame is of size $240 \time 240 \times 240$. They are stored as lateral slices of the tensor $\mathcal{X}_{\text{true}} \in \mathbb{R}^{240 \times 6 \times 240}$, and blurred by $\mathcal{A}\in \mathbb{R}^{240 \times 240 \times 240}$ in \eqref{za} with {\tt band} = 12 and $\sigma=2.5$ to obtain blur-contaminated, but noise-free frames stored as  lateral slices of $\mathcal{B}_{\text{true}} \in \mathbb{R}^{240 \times 6 \times 240}$. The blurred and noisy frames are generated by $\mathcal{B} = \mathcal{A}*\mathcal{X}_\text{true} + \mathcal{E}$, where $\mathcal{E}$ is defined above. The true fifth frame, and the blur- and noise-contaminated fifth frame are displayed in Figure \ref{Fig: 5s}. The restored fifth frames by the T-tEVD, T-tSVD, and RT-tSVD methods are shown in Figure \ref{Fig: 6s}.

\begin{figure}[htb]
\hspace{1cm}
\minipage{0.42\textwidth}
\includegraphics[width=\linewidth]{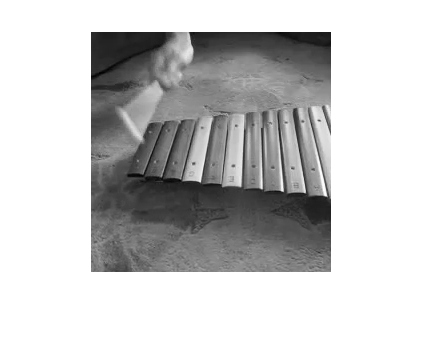} 
\endminipage\hfill \hspace{-2cm}
\minipage{0.42\textwidth}
\includegraphics[width=\linewidth]{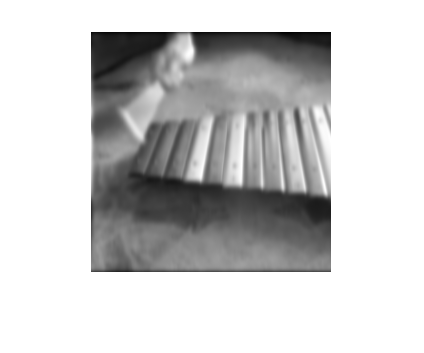}
\endminipage\hfill \hspace{-1cm} \vspace{-.9cm}
\caption{\small True fifth {\tt video} frame (left), and blurred and noisy fifth {\tt video} frame (right) for $\widetilde{\delta}=10^{-3}$.}
\label{Fig: 5s}
\end{figure}

\begin{figure}[!htb]
\hspace{-1cm}
\minipage{0.42\textwidth}
\includegraphics[width=\linewidth]{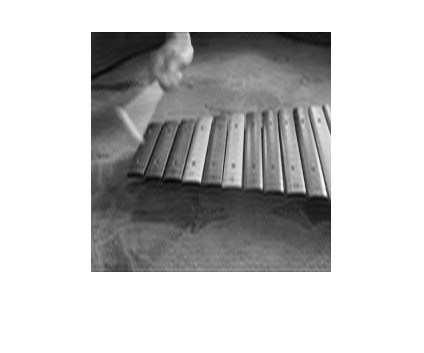} 
\endminipage\hfill \hspace{-2cm}
\minipage{0.42\textwidth}
\includegraphics[width=\linewidth]{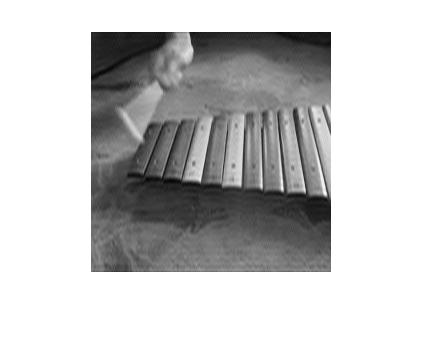}
\endminipage\hfill \hspace{-2cm}
\minipage{0.42\textwidth}
\includegraphics[width=\linewidth]{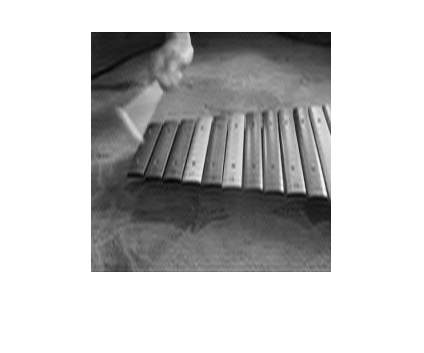}
\endminipage \hspace{-1cm} \vspace{-.9cm}
\caption{\small Reconstructed  fifth {\tt video} frame by the T-tEVD method after $72$ iterations (left), T-tSVD method after $72$ iterations (middle), and RT-tSVD method after $72$ iterations (right) for $\widetilde{\delta}=10^{-3}$.}
\label{Fig: 6s}
\end{figure}

Table \ref{Tab: 3s} shows that the {\tt nested}$\_$tGKT$_p$ method that is based on nested t-Krylov subspaces is competitive in terms of accuracy of the computed reconstruction and speed. This approach requires less than half of the computing time of the tGKT$_p$ method. The latter method requires the most CPU time and yields restoration of the highest quality for both noise levels. Among the truncated iteration based methods, the tGKB$_p$ method yields restoration of the best quality for $\widetilde{\delta}=10^{-3}$, while the T-tEVD and T-tSVD methods give reconstructions of the same quality, and of the highest quality for $\widetilde{\delta}=10^{-2}$. The t-Lanczos$_p$ method gives restoration of the worst quality and it is the fastest method for both noise levels. The quality of the reconstruction determined by this method is seen to improve when Tikhonov regularization is applied. The RT-tSVD method yields a reconstruction of higher quality than the T-tSVD method for $\widetilde{\delta}=10^{-3}$ and is faster than the latter method for both noise levels.

\begin{table}[h!]
\begin{center}
\begin{tabular}{cccccccc}
\cmidrule(lr){1-5}
\multicolumn{1}{c}{Noise level} & \multicolumn{1}{c}{Method}&\multicolumn{1}{c}{$k$} & \multicolumn{1}{c}{Relative error}& \multicolumn{1}{c}{CPU time (secs)}
\\ \cmidrule(lr){1-5}
\multirow{9}{4em}{$10^{-3}$} &tGKB$_p$ &(17,18,18,17,17,17)&$4.3850\cdot 10^{-2}$ &212.20 \\
&{\tt nested}$\_$tGKB$_p$ &22 & $4.4595\cdot 10^{-2}$ & 57.16\\
&t-Lanczos$_p$ &(5,5,5,5,5,5) & $6.1449\cdot 10^{-2}$ &8.59 \\
&T-tEVD&72& $4.5041\cdot 10^{-2}$ & 42.33 \\
&T-tSVD &72 & $4.5041\cdot 10^{-2}$ & 38.04 \\
&RT-tSVD&72 & $4.4962\cdot 10^{-2}$ & 21.67\\
&tGKT$_p$&(17,18,18,17,17,17)& $4.2660\cdot 10^{-2}$ & 238.71 \\
&{\tt nested}$\_$tGKT$_p$ &21& $4.3313\cdot 10^{-2}$ & 79.14\\
&t-LanczosTik$_p$&(5,5,5,5,5,5) & $4.7310\cdot 10^{-2}$ & 11.50 \\ \cmidrule(lr){1-5}
\multirow{9}{4em}{$10^{-2}$} &tGKB$_p$ &(4,4,4,4,4,4) &$7.6165\cdot 10^{-2}$ &12.27 \\
&{\tt nested}$\_$tGKB$_p$ &4 & $1.0741\cdot 10^{-1}$ &2.28 \\
&t-Lanczos$_p$ &(2,2,2,2,2,2) & $1.9410\cdot 10^{-1}$ &1.22\\
&T-tEVD&49& $7.3831\cdot 10^{-2}$ & 28.74 \\
&T-tSVD &49 & $7.3831\cdot 10^{-2}$ & 24.55 \\
&RT-tSVD&56 & $7.9567\cdot 10^{-2}$ & 16.36 \\
&tGKT$_p$&(4,4,4,4,4,4)& $6.4306\cdot 10^{-2}$ & 15.58 \\
&{\tt nested}$\_$tGKT$_p$ &4& $1.0690\cdot 10^{-1}$ & 7.62\\
&t-LanczosTik$_p$&(2,2,2,2,2,2) & $8.2925\cdot 10^{-2}$ & 3.07 \\ \cmidrule(lr){1-5}
\end{tabular}
\end{center} \vspace{-.5cm}
\caption{\small Results for tGKB$_p$, {\tt nested}$\_$tGKB$_p$, t-Lanczos$_p$, T-tEVD, T-tSVD, RT-tSVD, tGKT$_p$, {\tt nested}$\_$tGKT$_p$ and t-LanczosTik$_p$ methods when applied to the restoration of gray-scale {\tt video} frames.}
\label{Tab: 3s}
\end{table}

\end{Ex}

\section{Conclusions}\label{sec6}

This paper investigates the solution of linear discrete ill-posed problems for third order tensors by truncated iterations, where the truncation index that is determined by the discrepancy principle is the regularization parameter. A new randomized tensor singular value decomposition (R-tSVD) is proposed and applied to determine a low tubal rank tensor factorization. The performance of the R-tSVD depends on a prescribed error tolerance. The choice of the error tolerance depends on the level of noise in the data. This ensures that the sizes of the factors associated with the low tubal rank tensor factorization are not over-estimated. Randomization is used to speed up and improve the accuracy of the truncated tSVD method for solving tensor ill-posed problems. The randomized truncated tSVD (RT-tSVD) method is faster than the truncated tSVD (T-tSVD) in all examples considered, but its performance depends on the noise level. The truncated tensor eigenvalue decomposition (T-tEVD) method is seen to yield the same quality restoration and require the same number of iterations as the T-tSVD method. The latter method is faster but this may depend on the time spent sorting the eigentubes and corresponding eigenmatrices in order of descending magnitude.

The solution of a linear discrete ill-posed problem for third order tensors can be determined in terms of a few t-singular triplets and t-eigenpairs. The t-singular triplets are made up of singular tubes of the largest Frobenius norm and corresponding singular matrices, while the t-eigenpairs are the eigentubes of the largest Frobenius norm and associated eigenmatrices. The computation of a few of the largest t-singular triplets and t-eigenpairs by the t-product Golub-Khan bidiagonalization (tGKB) and t-product Lanczos (t-Lanczos) processes, respectively, can be inexpensive for t-linear discrete ill-posed problems. Computed examples suggest that knowledge of the t-singular triplets and t-eigenpairs may be quite informative when determining the solution of linear discrete ill-posed problems. Computing the t-singular triplets by a few steps of tGKB process is competitive when compared to the T-tSVD method. In general, the smaller the noise level, the more accurate the quality of restorations, and the more the t-singular triplets and t-eigenpairs required in order to satisfy the discrepancy principle.

The tGKB-type methods give the best or near-best quality restorations, while the t-Lanczos-type methods yield restorations of the worst quality for both noise levels. For $0.1\%$ noise, the tGKB$_p$ method is the slowest while the T-tEVD method is the slowest for $1 \%$ noise. Independently of the noise levels, the t-Lanczos-type methods are the fastest. Their performance as well as the performance of methods based on bidiagonalization process can be improved through the use of Tikhonov regularization. The quality of the computed approximate solutions can be marginally improved by truncating the iterations one step before the discrepancy principle is satisfied.

\end{document}